\RequirePackage{fix-cm}
\documentclass[coap]{svjour}                     
\smartqed  
\usepackage{graphicx}
\usepackage{mathtools}
\usepackage{enumitem}
\usepackage{bm}
\usepackage[ruled,linesnumbered,algo2e]{algorithm2e}
\SetKw{Kwin}{$\in$}
\SetKw{Kwpara}{in parallel}
\SetKw{KwOp}{Output:} 
\SetKwRepeat{Do}{do}{while}
\newlength\myindent
\makeatletter
\newcommand{\nosemic}{\renewcommand{\@endalgocfline}{\relax}}
\newcommand{\dosemic}{\renewcommand{\@endalgocfline}{\algocf@endline}}
\makeatother
\usepackage{wrapfig}
\usepackage[figuresright]{rotating}
\usepackage{amssymb}
\usepackage{hyperref}

\newcommand{\pluseq}{\mathrel{+}=}      
\usepackage[sort,numbers]{natbib}

\usepackage{nomencl}
\makenomenclature

\usepackage{multirow}
\usepackage{ragged2e}
\usepackage{graphicx}
\usepackage{textcomp}
\usepackage{lipsum}
\usepackage{xcolor}
\usepackage{array, booktabs, colortbl, tabularx}
\usepackage{pdflscape}
\usepackage{afterpage}
\usepackage{stfloats}
\usepackage{algpseudocode}
\usepackage{wrapfig}
\usepackage{subcaption}
\usepackage{caption}
\usepackage{array}
\newcolumntype{M}[1]{>{\centering\arraybackslash}m{#1}}

\DeclareMathOperator*{\argmin}{argmin}
\DeclareMathOperator*{\minimize}{minimize}
\DeclareMathOperator*{\maximize}{maximize}

\newcommand*\oline[1]{%
	\vbox{%
		\hrule height 0.5pt
		\kern0.35ex
		\hbox{%
			\kern-0.1em
			\ifmmode#1\else\ensuremath{#1}\fi
			\kern-0.1em
		}
	}
}

\makeatletter
\def\redefparbox{\def\@parboxrestore{\@arrayparboxrestore\let\\\@normalcr
		\if@minipage\expandafter\@gobbletwo\fi
		\@firstofone{\centering\casscparboxtest}}}
\def\casscparboxtest#1{%
	\ifx\rightskip#1\relax\expandafter\dimen@\else
	\expandafter\@secondoftwo
	\fi\@gobble{#1}}
\makeatother

\ExplSyntaxOn
\cs_gset_protected:Npn \__make_fig_caption:nn #1#2
{
	\skip_vertical:N \l_fig_abovecap_skip
	\parbox { \dim_eval:n { \l_fig_width_dim } }
	{
		\tl_use:N \l_fig_align_tl
		\sffamily \small \textbf{\color{scolor}#1:}~#2\par
	}
	\skip_vertical:N \l_fig_belowcap_skip
}
\cs_gset_protected:Npn \__make_tbl_caption:nn #1#2
{
	\skip_vertical:N \l_tbl_abovecap_skip
	\parbox{ \dim_eval:n { \l_tbl_width_dim } }
	{
		\tl_use:N \l_tbl_align_tl
		\sffamily \small \textbf{\color{scolor}#1}\par#2\par\vskip4pt
	}
	\skip_vertical:N \l_tbl_belowcap_skip
}
\ExplSyntaxOff

\usepackage{amsthm}
\journalname{XXXX Journal}
\begin{document}
\title{A Novel Column Generation Heuristic for Airline Crew Pairing Optimization with Large-scale Complex Flight Networks}
\author{Divyam Aggarwal\inst{1} \and Dhish Kumar Saxena\inst{1}\thanks{\emph{Corresponding author}; \emph{Email Address:} dhish.saxena@me.iitr.ac.in, dhishsaxena@gmail.com; \emph{Present Address:} Room No.-231, East Block, MIED, IIT Roorkee, Roorkee, Uttarakhand-247667, India; \emph{Phone:} +91-8218612326} \and Saaju Paulose\inst{2} \and Thomas B\"ack\inst{3} \and Michael Emmerich\inst{3}
}                     
%
%
\institute{Department of Mechanical \& Industrial Engineering (MIED), Indian Institute of Technology Roorkee, Roorkee, Uttarakhand-247667, India. \{daggarwal, dhish.saxena\}@me.iitr.ac.in \and
GE Digital Aviation Software, San Ramon, CA, USA. saaju.paulose@ge.com \and
Leiden Institute of Advanced Computer Science (LIACS), Leiden University, Niels Bohrweg 1, 2333 CA Leiden, the Netherlands. \{t.h.w.baeck, m.t.m.emmerich\}@liacs.leidenuniv.nl}
\date{}
%
\abstract{
    \textit{Crew Pairing Optimization} (CPO) is critical for an airlines’ business viability, given that the crew operating cost is second only to the fuel cost. CPO aims at generating a set of flight sequences (\textit{crew pairings}) to cover \textit{all} scheduled flights, at \textit{minimum} cost, while satisfying several \textit{legality} constraints. The state-of-the-art heavily relies on relaxing the underlying Integer Programming Problem into a Linear Programming Problem, that in turn is solved through the \textit{Column generation} (CG) technique. However, with the alarmingly expanding airlines’ operations, CPO is marred by the curse of dimensionality, rendering the \textit{exact} CG-implementations obsolete, and necessitating the \textit{heuristic} based CG-implementations. Yet, in literature, the much prevalent large-scale complex flight networks involving multiple – crew bases and/or hub-and-spoke sub-networks, largely remain uninvestigated. This paper proposes a novel \textit{CG heuristic}, which has enabled in-house development of an Airline Crew Pairing Optimizer (\textit{AirCROP}). The efficacy of the heuristic/\textit{AirCROP} has been tested on real-world, large-scale, complex network instances with over 4,200 flights, 15 crew bases and multiple hub-and-spoke sub-networks (resulting in a billion-plus possible pairings). Notably, this paper has a dedicated focus on the proposed CG heuristic (not the entire AirCROP framework) based on balancing \textit{random exploration of pairings}; \textit{exploitation of domain knowledge} (on optimal solution features); and \textit{utilization of the past computational \& search effort through archiving}. Though this paper has an airline context, the proposed CG heuristic may find wider applications across different domains, by serving as a template on how to utilize domain knowledge to better tackle combinatorial optimization problems.
} 
\maketitle
\section{Introduction} \label{intro}
\noindent Amongst major expenses of an airline, crew operating cost is the second-largest expense, after the fuel cost. For instance, 1.3 billion USD were spent on crew operations by American Airlines in 1991 \cite{gopalakrishnan2005airline}. For large airlines, even marginal improvements in crew operating cost may translate to savings worth millions of dollars annually. This has led to the recognition of Airline Crew Scheduling \cite{barnhart2003airline} as a critical planning activity with potential for huge cost-savings. Its aim is to generate crew schedules by assigning the crew members to a set of legal flight sequences, to cover a finite set of flights from its timetable while satisfying crew requirements of these flights. In the last three decades, it has received an unprecedented attention from the Operations Research (OR) society by way of comprehensive efforts to adapt the state-of-the-art optimization techniques to solve it. Conventionally, it is tackled sequentially through solutioning of two subproblems, including: (a) \textit{Crew Pairing Optimization Problem} (CPOP), where, the aim is to generate a set of flight sequences (each called a \textit{crew pairing}) to cover all the flights of an airlines' flight schedule, at minimum cost, while satisfying multiple \textit{legality} constraints linked to airline-specific regulations, federations' rules, labor laws, etc., and (b) \textit{Crew Rostering Problem} (or \textit{Crew Assignment Problem}), where, the aim is to assign crew members to optimally-derived crew pairings, while satisfying the corresponding crew needs. 
\par This research has contributed to in-house development of an Airline Crew Pairing Optimizer, named as \textit{AirCROP}\footnote{D. Aggarwal, D.K. Saxena, T. B\"ack, M. Emmerich, Crew Optimization, \textit{Netherlands Patent Application N2025010}, Feb. 2020.}, which has been tested and validated for real-world, large-scale and complex flight networks. However, this paper has a dedicated focus on the central building block for \textit{AirCROP} -- on \textit{how to intelligently explore the inestimably large search space of crew pairings (several billions possible), leading to minimization of cost}. The associated facts are highlighted below. CPOP belongs to the class of \textit{NP-hard}\footnote{For NP-hard (NP-complete) problems, no polynomial time algorithms on sequential computers are known up to now. However, verification of a solution could be accomplished efficiently, i.e., in polynomial time.} problems \cite{garey2002computers}, and its different facets include:
\begin{itemize}
	\item modeling: a CPOP is modeled either as a \textit{Set Covering Problem} (SCP), where, coverage of each flight in more than one pairing is permitted, or a \textit{Set Partitioning Problem} (SPP), where coverage of each flight is restricted to only one pairing.
	\item solution-architecture: for a crew pairing to be `operational' or `legal', it has to comply-with hundreds of \textit{legality} constraints (detailed in Section~\ref{sec:constraints}). For large-scale CPOPs, it is critically important to generate legal crew pairings in a timely manner. Based on either duty- or flight-network, several legal pairing generation approaches are available in the literature, which are reviewed in \citep{aggarwal2018large}. In this background, two CPOP solution architectures relate to how the legal pairing generation is invoked. One way is to completely enumerate all possible legal pairings \textit{a priori} CPOP-solutioning (finding a minimal-cost subset of those pairings covering all flights). This approach is mostly adopted for solving small-scale CPOPs, as generation and storage of all legal pairings for such problems is still computationally-tractable (see, e.g., \cite{beasley1996genetic, klabjan2001solving, zeren2012improved, deveci2018evolutionary}). Alternatively, the legal pairing generation is invoked over successive iterations of CPOP-solutioning, wherein, at each iteration, new legal pairings are enumerated to help CPOP solution improve as far as possible, before the next iteration is triggered (see, e.g., \cite{desaulniers1997crew, ahmadbeygi2009integer, duck2011implementing, muter2013solving, parmentier2020aircraft, desaulniers2020dynamic}). This approach is generally adopted for medium- to large-scale CPOPs, where millions/billions of legal pairings are possible, making their complete enumeration a priori CPOP-solutioning, computationally-intractable.
	\item solution-methodology: it is solved by use of either heuristic-based optimization techniques, or mathematical programming techniques, a brief review of which is presented below.
\end{itemize}
\subsection{Crew Pairing Optimization: Related Work} \label{sec:relatedwork}
As mentioned above, heuristic-based optimization techniques, and mathematical programming techniques are the two broad solutioning categories for the CPOPs. Notably, within heuristic-based techniques, the most widely adopted technique is Genetic Algorithms (GA), which are population-based, randomized-search heuristics, inspired by the theory of natural selection and genetics (details in \citep{goldberg2006genetic}). \citet{beasley1996genetic} is the first instance to customize GAs for solving a general class of SCPs, involving small-scale synthetic test cases (with just over $1,000$ rows and $10,000$ columns). The key facts around other GA-based CPOP solution approaches available, are summarized in Table~\ref{tab:GAlitRev}. 
\begin{table}[htbp]
	\small
	\centering \caption{An overview of GA-based CPOP solution approaches, proposed in the literature}
	\begin{center}
		\begin{tabular}{ccccc}
			\toprule
			\textbf{Literature Studies} & \textbf{Modeling} & \textbf{Timetable} & \textbf{Airline Test Cases*} & \textbf{Airlines}\\
			\midrule
			\citet{levine1996application} & Set Partitioning & - & 40R; 823F; 43,749P & - \\
			\citet{ozdemir2001flight} & Set Covering & Daily & 28R; 380F; 21,308P & Multiple Airlines \\
			\citet{kornilakis2002crew} & Set Covering & Monthly & 1R; 2,100F; 11,981P & Olympic Airways \\
			\citet{zeren2012improved} & Set Covering & Monthly & 1R; 710F; 3,308P & Turkish Airlines \\
			\citet{demirel2017novel} & Set Covering & Monthly & 6R; 1,002F; 1,121,408P & Turkish Airlines \\
			\citet{deveci2018evolutionary} & Set Covering & Monthly & 12R; 714F; 43,091P & Turkish Airlines \\
			\bottomrule
		\end{tabular}
		\\
		R represents the number of real-world test cases considered; F and P represents the maximum number of flights and pairings covered, therein. 
	\end{center}
	\label{tab:GAlitRev}
\end{table}
In general, these studies center around reasonably small-scale flight, for which relatively small number of pairings are possible. Though, \citet{kornilakis2002crew} and \citet{zeren2012improved} deal with $2100$ and $710$ flights respectively, they have considered only a subset of all possible legal pairings towards their reported solution. \citet{zeren2012improved} demonstrated that, while a highly customized-GA is efficient for solving small-scale CPOPs, it fails to solve large-scale CPOPs with the same search-efficiency. Furthermore, with reference to a $839$ flight-network with multiple hub-and-spoke sub-networks, \citet{aggarwal2020realworld} demonstrated that customized GAs are not so efficient in solving complex versions of even small-scale flight networks.
\par The use of mathematical programming techniques is linked to the fact that CPOPs modeled as SCPs/SPPs are inherently Integer Programming Problems (IPPs). However, it is found that the integer programs (IP) resulting from even small-scale CPOPs are so complex that it is computationally-impractical to solve them using standard IP-techniques \citep{anbil1991recent, kasirzadeh2017airline}. As an alternative, the integer constraint in CPOP/IPP is relaxed, leading to the corresponding Linear Programming Problem (LPP), solution to which is obtained by iteratively invoking an LPP solver. For the first iteration, the input to the LPP solver is any set of pairings which cover all the flights in the given schedule. For any subsequent iteration, the new input set comprises of the current LPP solution and a set of new pairings which by \textit{construct} promise improvement in the objective function. Such new pairings are generated using the \textit{Column Generation} (CG) technique, as part of the solution to a \textit{pricing sub-problem},  (\citet{desrosiers2005primer}). Once the optimal solution to the LPP is obtained, it is integerized towards a solution to the original CPOP/IPP. However, as cited by \cite{zeren2016novel}, the efficacy of the CG technique may be marred by-- tailing-off effect (slower convergence in the later LPP iterations), bang-bang effect (oscillation of dual variables from one extreme point to another), heading-in effect (poor dual information leading to generation of irrelevant columns in initial LPP iterations), etc. While, different CG stabilization techniques are available (\cite{du1999stabilized, lubbecke2005selected, lubbecke2010column}), the use of \textit{interior-point} methods (\cite{karmarkar1984new}) is gaining prominence. In terms of integerization, leading to the solution to the original CPOP, the literature points to prevalence of two approaches. One approach \citep{anbil1992global, zeren2016novel, parmentier2020aircraft} is to employ special performance boosting heuristics before applying a branch-and-bound algorithm (standard IP-technique). These heuristics help eliminate some pairings by exploiting their linear variables and by fixing specific flight-connections, or add some pairings before invoking the IP solver. The other approach is based on a \textit{branch-and-price} algorithm (originally proposed by \citet{barnhart1998branch}), in which CG is used to generate new legal pairings at each node of the branch-and-bound search tree \citep{vance1997heuristic, desaulniers2010airline, saddoune2013aircrew, quesnel2017new, quesnel2020branch, quesnel2020improving, desaulniers2020dynamic}.
\par For generating new pairings as part of the pricing sub-problem, the literature points to two types of strategies, namely, \textit{exact-} and \textit{heuristic-pricing} strategies. The aim of the exact-pricing strategy is to find only the pairings with most-negative reduced cost values, and if no pairing is found, then CG is terminated (see, e.g., \citet{vance1997airline, desaulniers2010airline, saddoune2013aircrew, desaulniers2020dynamic}). Here, the pricing sub-problem is modeled as a Shortest-Path Problem with Resource Constraints (SPPRC) \citep{irnich2005shortest}, and is solved using dynamic programming and labeling algorithms \citep{ahuja1993network}. \citet{minoux1984column} constructed the pricing sub-problem using a duty-tree, and applied a shortest path algorithm to find the most promising legal pairing (with most negative reduced cost value). In that, while finding the shortest path, i.e., the desired legal pairing, the cost of duties and their overnight-connections is reduced using the corresponding dual variables, and the pairing cost is assumed to be a linear sum of the these costs. This work is extended in \citet{lavoie1988new, barnhart1994column, desaulniers1997crew, vance1997heuristic}, where, primarily, the limitations linked to the linear construct of the pairing cost are addressed. \citet{saddoune2013aircrew} formulated individual pricing sub-problems for each combination of a day of planning horizon and a crew base, and defined each of them on a flight-based acyclic time-space network. \citet{desaulniers2020dynamic} adopted the same approach for their academic monthly instances (with over 7,527 flights \& 3 crew bases). However, for their industrial weekly (with over 8,711 flights \& 7 crew bases) and monthly (46,588 flights \& 7 crew bases) instances, the authors defined each pricing sub-problem on a duty-based acyclic time-space network, and solved them using a heuristic labeling algorithm with label-dominance rule \citep{irnich2005shortest}.
\par In the heuristic-pricing strategy, the aim is to find a subset of new pairings with negative reduced cost values either randomly or by exploiting problem-specific features \citep{marsten1994crew, zeren2016novel}. \citet{marsten1994crew} proposed an approach, where instead of solving the pricing sub-problem as SPPRC, the reduced cost of every possible legal pairing is explicitly checked. Notably, this approach is practically viable only in small- to medium-scale CPOPs, where the new pairings are reasonably sized. However, with the growing complexity of the flight networks over the last two decades, the utility of such CG implementations is impaired, paving the way for the domain-knowledge driven CG heuristics that focus upon a manageable, yet critical part of the overall pairings' search space. For instance, in an attempt to solve monthly flight schedules (with over $\approx$ $17,318$ flights \& 1 crew base) of Turkish Airlines, \citet{zeren2016novel} proposed a combination of both heuristic- and exact-pricing strategies, in which the former prioritizes the generation of a pairing set with considerably less number of deadhead, while the latter is solved as a shortest path IP problem. Notably, this IP formulation is not able to model all the pairing legality constraints, leading to the requirement of another heuristic for eliminating the generated infeasible pairings. Despite this progress, the applicability of heuristic-pricing strategy to the much prevalent and emergent \textit{complex flight networks}, characterized by \textit{multiple crew bases} and \textit{hub-and-spoke sub-networks}
, largely remains uninvestigated. This could plausibly be attributed to the fact that in such networks, the number of potential crew pairings grow exponentially with the number of flights. Such a research gap is all the more alarming, considering that the air traffic is expected to double up over the next $20$ years \citep{marisa2018airtravel}, where in, more and more airlines may rely on multiple crew bases and multiplicity of hubs. 
\subsection{Current Contributions}
This paper attempts to bridge the research gap through its proposition of a domain-knowledge driven CG heuristic capable of efficiently tackling large-scale and complex flight networks. In that, at any instant the proposed CG heuristic relies on balancing \textit{random exploration} (of pairings' space), \textit{exploitation of domain knowledge} (on optimal solution features) at a \textit{set} level and individual \textit{pairing} level, and utilization of past computational effort through  \textit{archiving} guided by the \textit{flight-pair level} information. The optimal solution features pursued here to guide the CG heuristic, relate to \textit{Deadhead reduction} (coverage of a particular flight in more than one pairing) and \textit{Crew Utilization enhancement} (the hours at work by a crew, out of the maximum permissible hours in a duty). The utility of the proposed CG heuristic is demonstrated on real-world, large-scale (over 4200 flights), complex flight network (over 15 crew bases and multiple hub-and-spoke sub-networks) of US-based client airlines of the research consortium's industrial partner. The proposed CG heuristic constitutes the kernel of a crew pairing optimization framework, named as \textit{AirCROP}, developed in-house as part of a larger research project, tested and validated by the industrial sponsor with reference to the best-practice solution for the underlying flight network data. Through experimental evidence using \textit{AirCROP}, it is established that the proposed CG heuristic based on integrated use of the four CG strategies, not only enables a lower-cost crew pairing solution in lesser time compared to any particular CG strategy or their sub-combinations, but also leads to a lower-cost solution compared to standard CG practice (with exact-pricing strategy).
\par The structure of the remaining paper is as follows. The formative concepts and the formulation of airline CPOP are presented in Section~\ref{sec:ACP}, followed by an overview of \textit{AirCROP} in Section~\ref{sec:aircrop}. The proposed CG heuristic is detailed in Section~\ref{sec:CGheuristic}. Subsequently, the experimental results are presented in Section~\ref{sec:exp}, while the paper concludes with Section~\ref{sec:conc}.
\section{Airline Crew Pairing Optimization Problem: Preliminaries and Formulation } \label{sec:ACP}
This section presents the prerequisites in terms of basic terminology, legality constraints and costing criteria for crew pairing, leading up to Airline CPOP formulation. 
\par
\begin{center}
	\begin{figure}[htbp]
		\centering{\includegraphics[width=0.62\columnwidth, keepaspectratio]{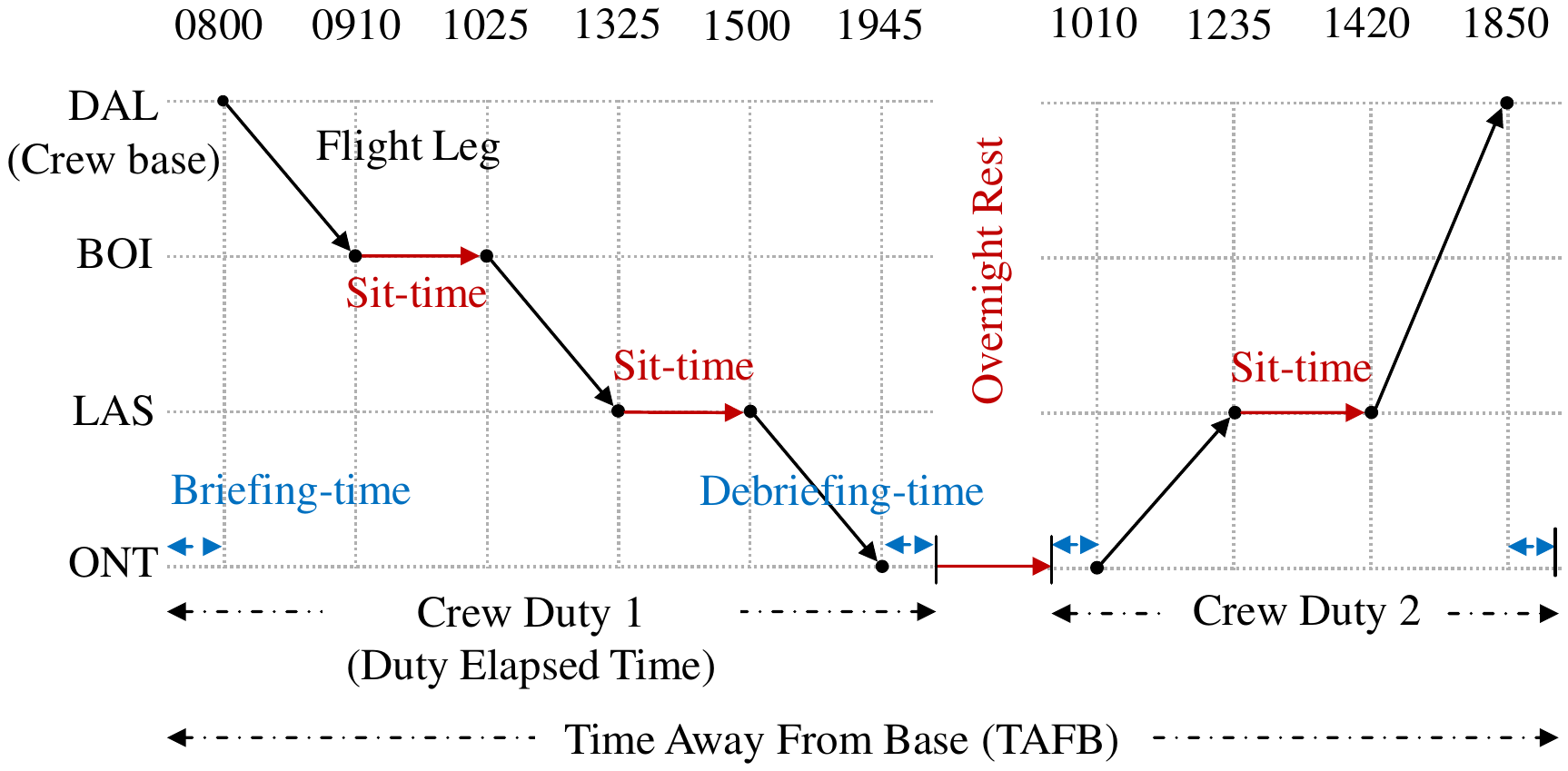}}
		\centering \caption{An example of a crew pairing starting from \textit{Dallas} (DAL) crew base.}
		\label{fig:pairing}
	\end{figure}
\end{center}

 A \textit{crew pairing} refers to a legal sequence of flights operated by a crew, that departs and arrives at a fixed (home) airport, called a \textit{crew base}. For instance, Figure~\ref{fig:pairing} illustrates a crew pairing with \textit{Dallas} (DAL) as the crew base.
The key terms used here, are discussed below. Within a crew pairing, the flight sequence covered in a single working day (not necessarily same as the calendar day) is called a \textit{crew duty}. In that, a small rest-period, provided between any two flights to facilitate aircraft changes by crew members, is called a \textit{sit-time} or a \textit{connection-time}. In contrast, \textit{overnight-rest} refers to a long rest-period after termination of the current duty and before commencement of the next crew duty. Two short-periods, provided in the beginning and ending of a crew duty, are called \textit{briefing} and \textit{de-briefing} time, respectively. The total time elapsed in a duty, including briefing and de-briefing times is called \textit{duty elapsed time}. Finally, the time elapsed since the start of a crew pairing is called the \textit{time away from base} (TAFB). At times, a crew is required to be transported to an airport to fly their next flight. In such situations, the crew is transported as passengers in another flight, flown by another crew. Such a flight is called a \textit{deadhead} or a \textit{deadhead flight} for the transported crew. The presence of deadhead flights affects an airline's profit in two ways. Firstly, the airline has to bear the loss of revenue on the passenger seats being occupied by the deadhead-ing crew. Secondly, the airline has to pay the hourly wages to the deadhead-ing crew even when they are not servicing the flight. In order to maximize the profits, airlines desire to minimize these deadheads as much as possible (ideally zero).
\subsection{Crew Pairing: Legality Constraints and Costing Criteria } \label{sec:constraints}
Multiple airline federations (Federal Aviation Administration, European Aviation Safety Agency, etc.) govern the safety of crew members, and, in turn, regulate the `legality' of a crew pairing. Along with this, several other legality constraints linked to airline-specific regulations, labor laws, etc. are required to be satisfied, for a crew pairing to be considered as legal. Broadly, these constraints could be classified as follows:
\begin{itemize}
	\item \textit{Connection-city constraint}: this requires the arrival airport of a flight to be the same as the departure airport of the next flight in the pairing sequence.
	\item \textit{Start- \& end-city constraint}: this requires the departure airport of the first flight and the arrival airport of the last flight, to be the same crew base.
	\item \textit{Sit-time \& overnight-rest constraint}: this requires the duration of the sit-times and overnight-rests to lie within the lower and upper bounds on the same, as set by airlines in accordance with the federations' regulation.
	\item \textit{Duty constraints}: these cater to the regulations on crew duty, including maximum limits on the permissible - number of duties in a pairing; number of flights in a crew duty; duty elapsed-time and its flying-time, etc.
	\item \textit{Special constraint}: Airlines desire to optimize their crew utilization, and formulate some special constraints, such as restricting a pairing which provides for overnights at the \textit{same city airports} (airports which are in the same city as the crew base), etc.
\end{itemize}
Given the multiplicity of these constraints, it is important to enable \textit{legal crew pairing generation} in a time-efficient manner, so that legal crew pairings are available as promptly as possible, when required during the optimization.
\par The cost of a crew pairing, in general, could be categorized into a \textit{flying cost} and a \textit{non-flying cost} (also called \textit{variable cost}). The former is the cost incurred in actually flying all given flights from an airline's schedule, and is calculated on an hourly basis. The latter is the cost incurred during the non-flying hours of a pairing, and could further be split into the following costs: (a) \textit{Hard cost} which includes the hotel cost, meal cost, and the  \textit{excess pay} -- the cost associated with the difference between the guaranteed hours of pay and the actual flying hours, and (b) \textit{Soft cost} which refers to the undesirable cost with associated aircraft changes (during flight-connections), etc.
\subsection{CPOP Formulation} \label{sec:model}
It has been mentioned earlier, that CPOP is modeled either as a set partitioning problem or as a set covering problem. The latter allows for more flexibility during its solutioning than the former, by accommodating deadhead flights. In that, for a given set of flights $\mathcal{F}$ (comprising of $F$ flights) that could be covered in multiple ways by a set of legal pairings $\mathcal{P}$ (comprising of $P$ pairings), the aim is to identify that subset of $\mathcal{P}$, say $\mathcal{P}_{IP}^*$, which not just covers \textit{each flight at least once}, but does it at a cost \textit{lower} than any alternative subset of the $\mathcal{P}$. The task of determining $\mathcal{P}_{IP}^* \subseteq \mathcal{P}$ is equivalent to each pairing $p_j \in \mathcal{P}$ being either included in $\mathcal{P}_{IP}^*$ (marked by the corresponding variable $x_j=1$) or excluded from  it (marked by the corresponding variable $x_j=0$). Notably, each $p_j$ is an $F$-dimensional vector, any element of which, say $a_{ij}$, is either $0$ or $1$. In that, if an $i^{th}$ flight ($f_i$) is covered by $p_j$, then $a_{ij}=1$, else $a_{ij}=0$. In this background, the CPOP  formulation (intrinsically an IPP), as employed in this paper, is presented below.
\begin{flalign}
	&\text{Min.}~ Z_{IP}  = \sum_{j=1}^{P} c_j x_j+ \psi_D \cdot \left(\sum_{i=1}^{F} \left(\sum_{j=1}^{P} a_{ij} x_j - 1\right) \right), \label{eq:obj} &\\
	&\text{subject to} \quad \sum_{j=1}^{P} a_{ij} x_{j} \geq 1,\quad ~~~~~\forall i \in \{1,2,...,F\} \label{eq:coverage}\\
	&\qquad \qquad \quad x_j \in \mathbb{Z}  = \{0,1\},~~~~\forall j \in \{1,2,...,P\} \label{eq:integrality}
\end{flalign}
\begin{flalign}
	\text{where},
	c_j &:~\text{the cost of a legal pairing }p_j,\nonumber &\\
	\psi_D &:~\text{an airline-defined penalty cost against each deadhead in the solution},\nonumber \\
	\quad a_{ij} &=~ 1,~\text{if flight}~f_i~\text{is covered in pairing}~p_j;~else~ 0 \nonumber \\
	x_j &=~ 1,~\text{if pairing}~p_j~\text{contributes to Minimum}~Z; ~else~ 0 \nonumber
\end{flalign}
In Equation~\ref{eq:obj}, the first cost component represents the sum of the individual cost of the pairings, while the second component represents the penalty cost for the deadheads in the solution (note, $(\sum_{j=1}^{P} a_{ij} x_j - 1)$ represents the number of deadheads corresponding to an $i^{th}$ flight in the solution). Notably, the above formulation assumes that the set of \textit{all} possible pairings, namely, $\mathcal{P}$, are available apriori, and the task is to determine the set $\mathcal{P}_{IP}^*$.  However, due to the practical challenge associated with generation of $\mathcal{P}$ apriori (as cited in Section~\ref{intro}), solution to the IPP is pursued in conjunction with the corresponding LPP (formulation deferred till Section~\ref{sec:CGheuristic}) assisted by CG technique.
\section{Overview of \textit{AirCROP} - the testbed for the Proposed CG heuristic} \label{sec:aircrop}
\begin{center}
	\begin{figure}[htbp]
		\centering{\includegraphics[width=0.8\columnwidth, keepaspectratio]{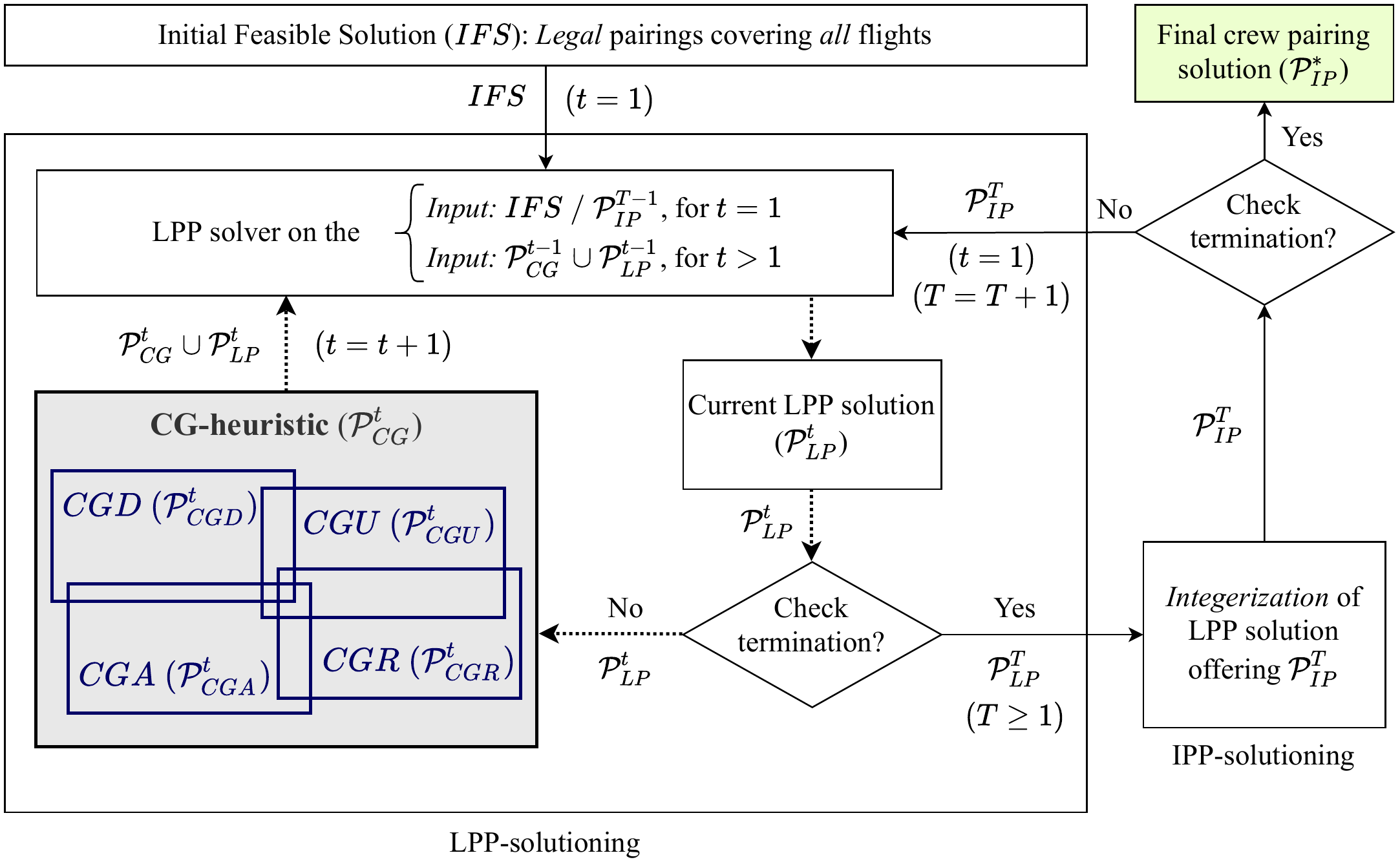}}
		\centering \caption{A schematic for \textit{AirCROP}, enabled by the proposed CG heuristic, comprising of the following strategies to generate fresh pairings $\mathcal{P}_{CG}^t$ at any LPP iteration $t$: \textit{Deadhead reduction} ($CGD$, generating $\mathcal{P}_{CGD}^t$), \textit{Crew Utilization enhancement} ($CGU$, generating $\mathcal{P}_{CGU}^t$), \textit{Archiving} ($CGA$, generating $\mathcal{P}_{CGA}^t$), and \textit{Random exploration} ($CGR$, generating $\mathcal{P}_{CGR}^t$). $T$ marks the counter to track the interactions between IPP- and LPP-solutioning.}
		\label{fig:aircrop}
	\end{figure}
\end{center}

As a prelude to the proposed CG heuristic (the core of this paper), an overview of \textit{AirCROP} is presented here, since it serves as a testbed for performance investigation of the heuristic. Its schematic in Figure~\ref{fig:aircrop}, highlights that it is developed by knitting the following modules -- \textit{Legal Crew Pairing Generation}, \textit{Initial Feasible Solution Generation}, and \textit{CG-driven LPP-solutioning}, interactively engaging with  \textit{IPP-solutioning}. The novelty in \textit{AirCROP} lies in not just how each module (briefed below) is conceived and implemented, but also in how these modules interact. An overview of these modules is presented in the following sub-sections. However, for a comprehensive discussion on their interactions, interested readers are referred to \citet{aggarwal2020aircrop}.
\subsection{Legal Crew Pairing Generation} \label{sec:LPGP}
This module helps generate pairings which satisfy all the \textit{legality} constraints. Such legal pairings enable generation of an initial feasible solution, and also feed real-time in to the proposed CG heuristic towards LPP-solutioning. Alarmingly, \citet{klabjan2001airline} have inferred based on empirical evidence, that typically 80\% of the CPOP-solutioning time may be consumed by legal pairings' generation. Considering that this challenge may only be augmented for complex flight networks, the following preemptive measures have been implemented, given which the total time spent on legal crew pairing generation ranges from $25$ to $50\%$ of the overall run-time for \textit{AirCROP}, depending on the scale and complexity of the underlying data set \citep{aggarwal2018large}:
\begin{itemize}
	\item \textit{adoption of duty-network:} here, all possible \textit{legal crew duties} are generated apriori by accounting for several duty constitutive constraints, including connection-city and sit-time constraints. Such \textit{duties} are then used to generate pairings, avoiding the need to re-evaluate several constraints, when pairings are needed in real-time.
	\item \textit{overnight-rest constraints' preprocessing:} here, a duty overnight-connection graph is pre-processed by satisfying the overnight-rest constraint between all possible pairs of apriori generated legal duties, eliminating the need for its repeated evaluation, when pairings are needed. 
	\item \textit{crew base driven parallelization:} recognizing, that duty constraints do vary across the different crew bases, legal pairing generation process is decomposed into independent sub-processes, capable of running in parallel on multiple cores of a single processing unit. 
\end{itemize}
\subsection{Initial Feasible Solution Generation} \label{sec:IFS}
This module aims to generate an \textit{initial feasible solution} (IFS) -- a manageable set of legal crew pairings covering all flights in a given schedule. For large-scale CPOPs, IFS generation standalone is computationally challenging as it constitutes an NP-complete problem. For \textit{AirCROP}, a time-efficient IFS-generation heuristic, namely, an \textit{Integer Programming based Divide-and-cover Heuristic} has been implemented \citep{aggarwal2020initializing}. It relies on: (i) a \textit{divide-and-cover} strategy to decompose the input flight schedule into smaller flight subsets, and (ii) use of \textit{integer programming} to find quality set of pairings for each decomposed flight subset.
\subsection{Interactive CG-driven LPP- and IPP-solutioning} \label{sec:optEng}
This is the core optimization module, which relies on intermittent interaction of submodules, as briefed below:
\begin{itemize}
	\item LPP-solutioning: as highlighted in Figure~\ref{fig:aircrop}, this phase entails several iterations ($t$) of an LPP solver until the termination criterion is not met. In that, if the cost improvements over a pre-specified number of successive LPP iterations falls below a pre-specified threshold (highlighted in Section~\ref{sec:computational-setup}), this phase is terminated. In the initial iteration ($t=1$), the input set is nothing but the IFS. In any subsequent iteration $t$, the input set comprises of the previous iteration's-- LPP solution $\mathcal{P}_{LP}^{t-1}$ and a fresh set of pairings generated/re-inducted using the proposed CG heuristic ($\mathcal{P}_{CG}^{t-1}$). This input pairing set is fed into the LPP solver to find the current $\mathcal{P}_{LP}^t$, which is further used by the proposed CG heuristic to generate a fresh set of pairings $\mathcal{P}_{CG}^{t}$. While the details of the proposed CG heuristic are presented in the next section, it may be noted here, that the composition of $\mathcal{P}_{CG}^t$ is guided by a two-pronged strategy of exploring the \textit{new} pairings' space on the one hand, and utilizing the efficient pairings encountered over the \textit{past} LPP iterations on the other. In that: 
	\begin{itemize}
		\item exploring the \textit{new} pairings' space: this is guided by three strategies, namely, \textit{deadhead reduction} ($CGD$), \textit{crew utilization enhancement} ($CGU$), and \textit{random exploration} of pairings' space ($CGR$). Notably, the former two strategies pursue the desired optimal solution features. 
		\item utilizing the efficient pairings from the \textit{past}: this is guided by the \textit{archiving} strategy ($CGA$) which utilizes the dual variables of the current LPP solution to identify the promising flight-pairs, corresponding to which efficient pairings from the archive of the past LPP iterations, are extracted.
	\end{itemize}
	\item IPP-solutioning: once a particular phase of LPP-solutioning terminates, the LPP solution is passed over for its integerization, towards which the branch-and-bound algorithm is used. Being a minimization problem, this algorithm maintains the valid lower bound (LPP solution) and the best upper bound (IPP solution) at each node of its search-tree, and it terminates itself if the difference between these two bounds (\textit{MIP gap)} becomes zero. However, the branch-and-bound search on a large-scale IPP is time-intensive. Hence, the algorithm is terminated at a pre-defined time limit (highlighted in Section~\ref{sec:computational-setup}), if it does not self-terminate apriori. Subsequently, the corresponding IPP solution is passed back as the input for the next phase of LPP-solutioning.
\end{itemize}
Notably, the LPP- and IPP-solutioning interactions are tracked through the $T$ count, and in principle \textit{AirCROP} is terminated when the cost of the IPP solution ($Z_{IP}^T$ corresponding to  $\mathcal{P}_{IP}^T$) matches the cost of its input LPP solution ($Z_{LP}^T$ corresponding to $\mathcal{P}_{LP}^T$). However, for practical considerations on the time limit, if the IPP and LPP costs do not conform with each other despite a pre-specified number of LPP-IPP interactions, or up to a pre-specified limit on the total run-time, the termination of \textit{AirCROP} is enforced (these pre-specified numbers are highlighted in Section~\ref{sec:computational-setup}). It may also be noted that aligned with the scope of this paper, all features and associated parameters in \textit{AirCROP} but for the combinations of CG strategies, are uniformly maintained for experimentation and results. 
\section{Domain-knowledge Inspired Column Generation Heuristic} \label{sec:CGheuristic}
The aim in this paper is to demonstrate the efficacy of the proposed CG heuristic towards fulfilling the overarching aim of solving the CPOP/IPP problem (Equations~\ref{eq:obj} to \ref{eq:integrality}) within the \textit{AirCROP} architecture (Figure~\ref{fig:aircrop}). Notably, the prerequisite for solving the IPP formulation is that the set of all possible legal pairings $\mathcal{P}$ is generated apriori. Considering, this may be computationally-intractable for large-scale/complex flight networks, \textit{AirCROP} relies on first solving the corresponding LPP (formulated below), and integerization of the resulting LPP solution, over multiple such interactions. Solution to LPP is obtained by iteratively invoking an LPP solver. In the first iteration, the input to an LPP solver is any set of pairings covering all the scheduled flights. In any subsequent iteration, the new input set comprises of the previous iteration's-- LPP solution and a set of fresh pairings promising further improvement in the objective function, as generated/re-inducted using the proposed CG heuristic. In other words, \textit{the quality of the final LPP solution is directly dependent on the efficacy of the CG heuristic}. Furthermore, whenever an LPP solution is passed on for its \textit{integerization} (IPP-solutioning), the objective function value for the LPP solution marks the lower bound for the corresponding Integer solution. In other words, an integer solution can not have a better objective function value than the underlying LPP solution. It implies that after several interactions of LPP- and IPP-solutioning, \textit{the quality of the final IPP solution would also depend on the search efficiency in the LPP solutioning phase, at the core of which lies the CG heuristic. This is a testament to the importance of the proposed CG heuristic and this paper}.
\par The details and associated nitty-gritty of the proposed CG heuristic can be best explained at the level of each LPP iteration ($t$), that could be perceived as being implemented in three phases\footnote{For ease of reference, the notations introduced here are kept independent of the iteration counter $t$. However, for subsequent pseudocodes and corresponding discussions with reference to a particular LPP iteration, these notations are superscripted by $t$}. In that:
\begin{itemize}
	\item firstly, for a given set of input pairings, LPP solver is invoked on the \textit{primal} form of the LPP formulation (Equations~\ref{eq:obj2}~to~\ref{eq:integrality1}), and the LPP solution is obtained, in that, each pairing $p_j$ in the input set gets assigned a value $x_j$ varying from $0$ to $1$. All the $x_j$s together constitute the \textit{primal} variable vector, notated as $X=[x_1~x_2~x_3~...~x_P]^{\mathsf{T}}$. Furthermore, the elements in $X$ with non-zero values ($x_j \neq 0$), and the set of underlying pairings are notated as $X_{LP}$ and $\mathcal{P}_{LP}$, respectively.   
	\begin{flalign}
		&\text{Min.}~ Z_{LP}^p = \sum_{j=1}^{P} c_j x_j+ \psi_D \cdot \left(\sum_{i=1}^{F} \left(\sum_{j=1}^{P} a_{ij} x_j - 1\right) \right)=\sum_{j=1}^{P} \left(c_j + \psi_D \cdot \sum_{i=1}^{F} a_{ij} \right) x_j - F \cdot \psi_D, &\label{eq:obj2} \\
		&\text{subject to} \quad \sum_{j=1}^{P} a_{ij} x_{j} \geq 1,\qquad \quad \forall i \in \{1,2,...,F\} \label{eq:coverage2}\\
		&\qquad \qquad \quad x_j \in \mathbb{R} = [0,1],\qquad \forall j \in \{1,2,...,P\} \label{eq:integrality1}
	\end{flalign}
	Notably, the contribution of each pairing to the LPP solution (Equation~\ref{eq:integrality1}), could be effectively treated as:
	\begin{flalign}
		&\qquad \qquad \quad x_j \in \mathbb{R} \geq 0,\qquad \quad~~ \forall j \in \{1,2,...,P\} &\label{eq:integrality2}
	\end{flalign}
	The rationale is that minimization of $Z_{LP}^p$ will always lead to a solution with all variables $x_j \leq 1$ \citep{vazirani2013approximation}. Hence, the primal form of the LPP is represented by Equations~\ref{eq:obj2},~\ref{eq:coverage2}~\&~\ref{eq:integrality2}. Here, the use of Equation~\ref{eq:integrality2} instead of Equation~\ref{eq:integrality1}, helps in getting rid of an additional dual variable in the dual form of the LPP formulation, as defined in the subsequent text.
	\item secondly, the pairing set $\mathcal{P}_{LP}$ is fed to the \textit{dual} form of the LPP formulation (Equations~\ref{eq:dualobj}~to~\ref{eq:dualvariables}), and LPP solver is invoked to yield the \textit{dual} variable vector $Y=[y_1~y_2~y_3~...~y_F]^{\mathsf{T}}$ corresponding to the flight constraints. 
	\begin{flalign}
		&\text{Max.}~ Z_{LP}^d = \sum_{i=1}^{F} y_i - F \cdot \psi_D, \label{eq:dualobj} &\\
		&\text{subject to} \quad \sum_{i=1}^{F} a_{ij} y_i \leq \left(c_j + \psi_D \cdot \sum_{i=1}^{F} a_{ij}\right),~~~~\forall j \in \{1,2,...,P_{LP}\} \label{eq:dualconstraints}\\
		& \qquad \qquad \qquad \quad~~y_i \in \mathbb{R} \geq 0,\qquad \qquad \qquad~~~~\forall i \in \{1,2,...,F\} \label{eq:dualvariables}\\
		&\text{where}, \quad P_{LP} :~\text{is the number of pairings in the set}~  \mathcal{P}_{LP} &  \nonumber \\
		&\qquad \qquad y_i :~\text{dual variable, corresponding to an $i^{th}$ flight-coverage constraint}, & \nonumber 
	\end{flalign}
	%
	\item availability of the \textit{dual} variables paves way for the \textit{pricing sub-problem}, where the aim is to generate those legal pairings which if included as part of the input to the next LPP iteration, promise a better $Z_{LP}^p$ than the current. The \textit{standard} CG technique (with exact-pricing strategy, as mentioned in Section~\ref{sec:relatedwork}) identifies only the pairings with most negative \textit{reduced cost} values, where, the reduced cost ($\mu_j$) is, as defined in Equation~\ref{eq:redCost}.
    \begin{flalign}
    	&\mu_j = \left(c_j + \psi_D \cdot \sum_{i=1}^{F} a_{ij}\right) - \mu d_j,~\text{where,}~\mu d_j = \sum_{i=1}^{F} a_{ij} \cdot y_i~~(~\text{represents dual cost component of}~\mu_j)&\label{eq:redCost}
    \end{flalign}
	Here, the potential of a pairing to further reduce $Z_{LP}^p$ is in proportion to the magnitude of its $\mu_j$ value. Notably, \textit{the distinctive contribution of the proposed CG heuristic} lies in the fact, that from within the larger pool of pairings with negative $\mu_j$, besides selecting pairings randomly, it also selects pairings in a guided manner. In that, selection of such pairings is guided by optimal solution features at a \textit{set level} and individual \textit{pairing level}, and utilization of the past computational effort. Once a set of such promising pairings meets a pre-defined size, it is merged with the existing $\mathcal{P}_{LP}$, and fed as the input for the next LPP iteration. When repeated, this procedure leads to a near-optimal LPP solution which could be integerized.
\end{itemize}   
In essence, while the \textit{standard} CG technique prioritizes the generation of only the pairings with most negative $\mu_j$ to qualify as an input for the next LPP iteration, the proposed CG heuristic prioritizes the generation of some pairings over the others amongst the pool of pairings with negative $\mu_j$, through its multi-pronged strategies, including:
\begin{itemize}
	\item Deadhead reduction strategy, $CGD$: prioritizes a set of legal pairings that is characterized by \textit{low deadheads} - a feature which domain knowledge recommends for optimality at a \textit{set level}.
	\item Crew Utilization enhancement strategy, $CGU$: prioritizes a set of legal pairings, each member of which is characterized by \textit{high crew utilization} - a feature which domain knowledge recommends for optimality at a \textit{pairing level}.
	\item Archiving strategy, $CGA$: prioritizes a set of legal pairings comprising of those flight-pairs, which as per the existing LPP solution, bear better potential for improvement in the objective function. Such a set, originating from the \textit{flight-pair level} information, is extracted from an archive of previously generated pairings. In doing so, this strategy facilitates utilization of the past computational efforts, by providing an opportunity for an earlier generated legal pairing to be re-inducted in the current pairing pool. 
	\item Random exploration strategy, $CGR$: unlike $CGU$ and $CGD$ which are guided by the desired optimal solution features, and $CGA$ which is guided by the existing LPP solution features, this strategy pursues random and unbiased exploration of pairings' space.  
\end{itemize}
In this background, the details of each strategy leading upto the proposed CG heuristic are presented below.
\subsection{Deadhead Reduction Strategy (\texorpdfstring{$CGD$}{CGD})} \label{sec:CGD}
As mentioned earlier, a deadhead flight adversely affects an airline's profit. Hence, airlines desire to minimize the deadhead flights as much as possible (ideally zero). Towards it, while defining the costing criteria for a pairing, airlines include a high penalty cost for each deadhead flight present in a crew pairing solution. It implies that while pursuing minimization of the total cost, the resulting solution shall also be characterized by minimal possible deadheads. To exploit this optimality feature, the \textit{deadhead reduction} strategy ($CGD$) proposed here, aims to generate a new pairing set, $\mathcal{P}^{t}_{CGD}$, which provides an alternative way to cover the flights involved in a subset of $\mathcal{P}^{t}_{LP}$, while also ensuring that some of these flights get covered with zero deadheads. This strategy promises propogation of the zero-deadhead feature over successive LPP iterations for the following reasons:
\begin{itemize}
	\item $\mathcal{P}^{t}_{CGD}$ alongside $\mathcal{P}^{t}_{LP}$ forms a part of the input for the next LPP iteration
	\item $\mathcal{P}^{t}_{CGD}$ provides a scope for better coverage (zero-deadhead) of some flights, compared to $\mathcal{P}^{t}_{LP}$
	\item $\mathcal{P}^{t}_{CGD}$ may focus on zero-deadhead coverage for different flights in different LPP iterations
\end{itemize}
\par
\begin{figure}[htbp]
	\centering
	\includegraphics[width=0.75\linewidth, keepaspectratio]{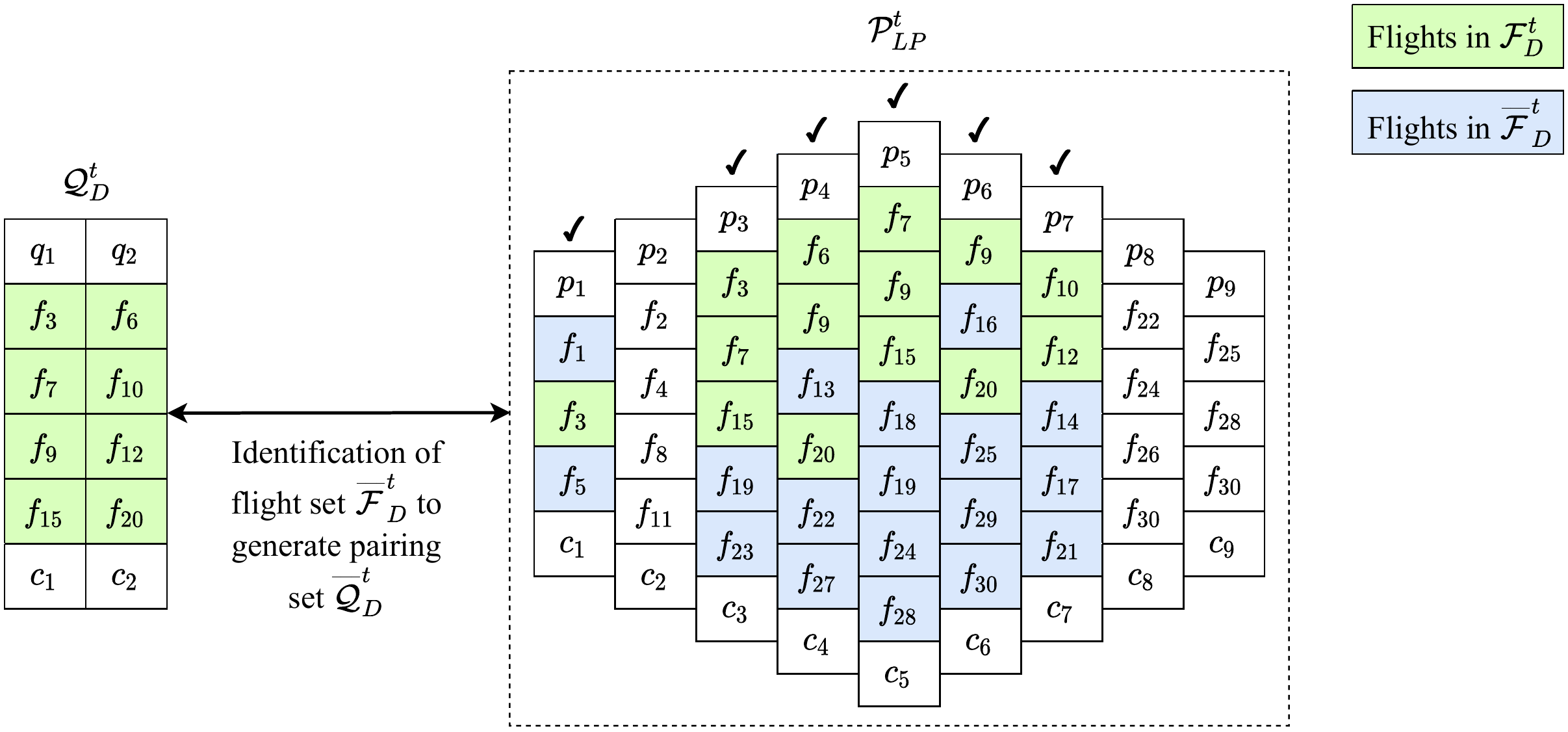}
	\caption{A symbolic depiction of how the deadhead reduction strategy provides for an alternative way to cover the flights involved in a subset of $\mathcal{P}^{t}_{LP}$, while also ensuring that some of these flights get covered with zero deadheads}.
	\label{fig:CGD}
\end{figure}

 The implementation of this strategy at iteration $t$, formalized in Algorithm~\ref{alg:CGD}, has been explained below in conjunction with Figure~\ref{fig:CGD} (where, a sample $\mathcal{P}^{t}_{LP}$ comprises of pairings $p_{1}$ to $p_{9}$, to cover the scheduled flights $f_1$ to $f_{30}$).
\begin{algorithm2e}[htbp]
	\small
	\DontPrintSemicolon
	\SetKwComment{Comment}{$\triangleright$\ }{}
	\KwIn{$\mathcal{D}_{all},~\mathcal{P}^{t}_{LP},~Y^{t},~Th_D$}
	\KwOut{$\mathcal{P}^{t}_{CGD}$}
	$Random \gets$ Generate a random integer between 0 \& $Th_D$\;
	$\mathcal{D}_{D}^t \gets$ Select $Random$ number of duties, with-respect-to each crew base, from $\mathcal{D}_{all}$ without replacement\;
	Generate all possible pairings using $\mathcal{D}_{D}^t$ and compute their $\mu^{t}_{j}$\Comment*{Use Equation~\ref{eq:redCost}}
	Select only the pairings with negative $\mu^{t}_j$\;
	Sort pairings in ascending order with-respect-to their $\mu^{t}_j$\;
	$\mathcal{Q}_{D}^{t} \gets$ Select a zero-deadhead pairing set from the sorted pairings\;
	$\mathcal{F}_{D}^{t} \gets$ Identify flights covered by pairings in $\mathcal{Q}_{D}^{t}$\;
	\For{each pairing $p_j \in \mathcal{P}^{t}_{LP}$}{
		\If{pairing $p_j$ covers one or more flights from $\mathcal{F}_{D}^{t}$}{
			$\hspace{1.5pt}\oline{\mathcal{F}}_{D}^{~t} \gets$ Add only those flights from $p_i$ that are not present in $\mathcal{F}_{D}^{t}$\;
		}
	}
	Generate all possible pairings using $\hspace{1.5pt}\oline{\mathcal{F}}_{D}^{~t}$ and compute their $\mu^{t}_{j}$\Comment*{Use Equation~\ref{eq:redCost}}
	$\hspace{1.5pt}\oline{\mathcal{Q}}_{D}^{~t} \gets$ Select only the pairings with negative $\mu^{t}_j$\;
	$\mathcal{P}^{t}_{CGD} \gets \mathcal{Q}_{D}^{t} \cup \hspace{1.5pt}\oline{\mathcal{Q}}_{D}^{~t}$\;
	\textbf{return} $\mathcal{P}^{t}_{CGD}$\;
	\caption{Procedure for the $CGD$ strategy, at an LPP iteration $t$}
	\label{alg:CGD}
\end{algorithm2e}
The initial step of $CGD$ involves random selection of a subset of legal duties, say $\mathcal{D}^t_D$, from the set of all legal duties, say $\mathcal{D}_{all}$ (line 2). The maximum number of duties contributing to $\mathcal{D}^t_D$ with-respect-to each crew base, are limited by a pre-defined threshold, say $Th_D$ (the criterion for its setting is highlighted in Section~\ref{sec:parameter-setting}). The actual number of duties contributing to $\mathcal{D}^t_D$ with-respect-to each crew base, are guided by $Random$ -- a random integer between 0 \& $Th_D$ (line 1). Subsequently, all possible pairings corresponding to $\mathcal{D}^t_D$ are generated, of which only those with negative reduced cost ($\mu_j^t$) are selected (lines 3 \& 4). These pairings are then sorted in ascending order with respect to their $\mu_j^t$ values (line 5), and the largest sized subset characterized by zero deadheads, say $\mathcal{Q}_D^t$, is extracted using a deterministic approach (line 6). In that, starting with the pairing having the most-negative $\mu_j^t$ (maximum potential for improvement in the objective), iterations are performed over the sorted pairings, and only those pairings are selected which do not cover an already covered flight. The above sorting is important, since the composition of the zero-deadhead subset is dependent on the order in which pairings are evaluated. Let the set of flights covered by pairings in $\mathcal{Q}_D^t$ be denoted by $\mathcal{F}_D^t$ (line 7). Subsequently, $\mathcal{P}_{LP}^t$ is scrolled through (lines 8 to 12), and only the pairings covering one or more flights in $\mathcal{F}^t_{D}$ are flagged. This enables the composition of $\hspace{1.5pt}\oline{\mathcal{F}}^{~t}_{D}$ (line 10) as the set of complementary flights -- those which are not present in $\mathcal{F}^t_{D}$. For instance, in Figure~\ref{fig:CGD}, $\mathcal{Q}^t_{D}=\{q_1,q_2\}$, $\mathcal{F}^t_{D}=\{f_3,~f_6,~f_7,~f_9,~f_{10},~f_{12},~f_{15},~f_{20}\}$ (flights in the green colour); pairings $p_{1},~p_{3}~\ldots~p_{7}$ are flagged; and $\hspace{1.5pt}\oline{\mathcal{F}}^{~t}_{D}= \{f_1,~f_5,~f_{13},~f_{14},~f_{16}~\ldots ~f_{19},~f_{21}~\ldots~f_{25},~f_{27}~\ldots~f_{30}\}$ (flights in the blue colour). Subsequently, all possible pairings corresponding to $\hspace{1.5pt}\oline{\mathcal{F}}^{~t}_{D}$ are generated, of which only those having negative $\mu_j^t$ values, constitute the set  $\hspace{1.5pt}\oline{\mathcal{Q}}^{~t}_{D}$ (lines 13 \& 14). Finally, the pairing sets  $\mathcal{Q}^t_{D}~\&~\hspace{1.5pt}\oline{\mathcal{Q}}^{~t}_{D}$ jointly constitute $\mathcal{P}^{t}_{CGD}$ (line 15). 
\subsection{Crew Utilization Enhancement Strategy (\texorpdfstring{$CGU$}{CGU})} \label{sec:CGU}
The degree of crew utilization in a pairing can be assessed by the number of hours a crew is at work, out of the maximum permissible working hours, in each of its constituent crew duties. Notably, higher the crew working hours in each duty of a pairing, lower will be its excess pay (Section~\ref{sec:constraints}). Hence, minimization of the total crew cost is directly linked to the presence of pairings with high crew utilization in the final LPP solution. To exploit this optimality feature, the \textit{crew utilization enhancement} strategy ($CGU$) proposed here:
\begin{itemize}
	\item first introduces a new measure, namely, \textit{crew utilization ratio} ($\gamma_j$, Equation~\ref{eq:crewUtilityRatio}) to quantify the degree of crew utilization in a pairing $p_j$ at any instant \begin{flalign}
		\gamma_j = \frac{1}{\text{Number of duties in }p_j} \cdot \sum_{d \in p_j} \frac{\text{Working hours in duty} ~d}{\text{Permissible hours of duty } d}  \label{eq:crewUtilityRatio}
	\end{flalign}
	\item identifies with reference to $\mathcal{P}^{t}_{LP}$, pairings characterized by high dual cost component ($\mu d_j^t$, Equation~\ref{eq:redCost}), reflecting in turn on the constitutive flights with high value of dual variables ($y_i^t$, Equation~\ref{eq:redCost})), and hence, on the potential of these flights to generate new pairings with more negative $\mu_j^t$
	\item utilizes the flight information to generate promising pairings, and picks the ones with high $\gamma_j^t$ to constitute $\mathcal{P}^{t}_{CGU}$. 
\end{itemize}
This strategy promises propogation of the higher crew utilization ratio over successive LPP iterations, given that in each LPP iteration $\mathcal{P}^{t}_{CGU}$ alongside $\mathcal{P}^{t}_{LP}$ forms a part of the input for the next iteration.
\par
\begin{figure}[htbp]
	\centering
	\includegraphics[width=0.75\linewidth, keepaspectratio]{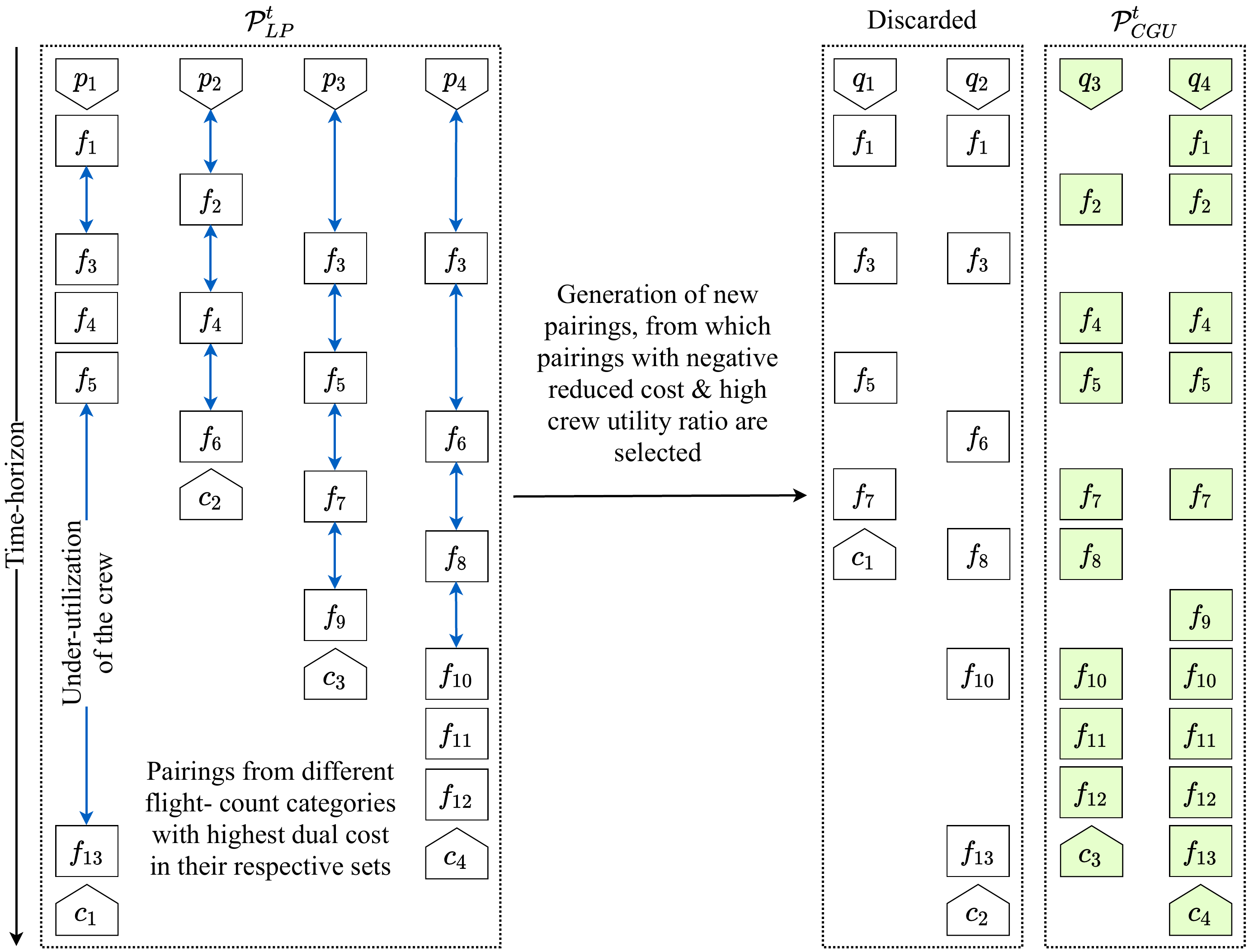}
	\caption{A symbolic depiction of how the crew utilization enhancement strategy provides for novel flight-connections, enabling generation of new pairings with high crew utilization ratio}
	\label{fig:CGU}
\end{figure}

 The implementation of this strategy at iteration $t$, formalized in Algorithm~\ref{alg:CGU}, has been explained below in conjunction with Figure~\ref{fig:CGU} (where, a sample $\mathcal{P}^{t}_{LP}$ comprises of pairings $p_{1}$ to $p_{4}$, covering the flights $f_1$ to $f_{13}$). 
\begin{algorithm2e}[htbp]
	\small
	\DontPrintSemicolon
	\SetKwComment{Comment}{$\triangleright$\ }{}
	\KwIn{$\mathcal{P}^{t}_{LP},~X_{LP}^t,~Y^{t},~Th_U$}
	\KwOut{$\mathcal{P}^{t}_{CGU}$}
	\nosemic 
	$Dictionary \gets$ Categorize pairings in $\mathcal{P}^{t}_{LP}$ with-respect-to flight-counts\;
	$Random \gets$ Generate a random integer between 0 \& $Th_U$\;
	\While{$|\mathcal{F}_U^{t}| < Random$}{
		Select a random $(k,\mathcal{P}_k^t)$ pair from $Dictionary$ without replacement\;
		Compute $\mu d_j^{t}$ values for pairings in $\mathcal{P}^{t}_{k}$\Comment*{Use Equation~\ref{eq:redCost}}
		Sort pairings in $\mathcal{P}^{t}_{k}$ in descending order with-respect-to their $\mu d^{t}_j$ \& $x_j~(\in X_{LP}^t)$\;
		$Counter \gets 0$\;
		\For{each pairing $p_j \in \mathcal{P}^{t}_{k}$}{
			$\mathcal{F}_U^{t} \gets$ Add flights covered in $p_j$\;
			$Counter \pluseq$ Number of flights in $p_j$\;
			\If{$Counter \geq Random /(k_{max}-k_{min})$}{
				Break\;
			}
		}
	}
	Generate all possible pairings using $\mathcal{F}_U^{t}$ and select only the pairings with negative $\mu^t_j$\;
	Compute $\gamma_j^t$ of the selected pairings and find the median value \Comment*{Use Equation~\ref{eq:crewUtilityRatio}}
	$\mathcal{P}^{t}_{CGU} \gets$ Select only the pairings with $\gamma_j^t \geq$ the median value\;
	\textbf{return} $\mathcal{P}^{t}_{CGU}$\;
	\caption{Procedure for the $CGU$ strategy, at an LPP iteration $t$}
	\label{alg:CGU}
\end{algorithm2e}
In that:
\begin{itemize} 
	\item firstly, all pairings in $\mathcal{P}^{t}_{LP}$ are categorized with-respect-to the number of flights covered in them (flight-count), denoted by $k$. The range of $k$ is given by the upper and lower bounds, denoted by $k_{max}$ and $k_{min}$, respectively
	\item secondly, a $Dictionary$ is constituted, each element of which can be given by $(k,\mathcal{P}_k^t)$, where, $\mathcal{P}_k^t$ represents the set of pairings corresponding to flight count $k$ (line 1). The rationale is to identify promising flights corresponding to each $k$, so that their collective set ($\mathcal{F}_U^t$, formed over all the $k$s), has a better chance to offer pairings with higher flight-count and crew utilization ratio. The size of $\mathcal{F}_U^t$ is controlled by $Random$ -- a random integer between 0 and a pre-defined upper threshold, namely, $Th_U$ (the criterion for its setting is highlighted in Section~\ref{sec:parameter-setting}).
	\item subsequently, $\mathcal{F}_U^t$ is composed iteratively, where each iteration involves:
	\begin{itemize}
		\item random selection of a $(k,\mathcal{P}_k^t)$ category from the $Dictionary$ without replacement (line 4)
		\item computation of each pairing's $\mu d_j^t$ (line 5), and subsequently their sorting by descending order of their $\mu d_j^t$, or by their primal variables $x_j^t$ in case of a tie in $\mu d_j^t$ values (line 6)
		\item sequential scrolling of the sorted pairings, and assignment of each constitutive flight to $\mathcal{F}_U^t$, until the flights add up to the limit permissible for each $k$ category  (lines 7 to 14). 
	\end{itemize}
	\item when the size of $\mathcal{F}_U^t$ exceeds its permissible limit, the iterations are terminated (line 15). Then, all possible pairings with the flights in $\mathcal{F}_U^t$ are generated, and only those with negative $\mu_j^t$s are picked, for computation of $\gamma_j^t$s (lines 16 \& 17). Finally, the pairings with above-median $\gamma_j^t$ values constitute the desired new pairing set, namely $\mathcal{P}_{CGU}^t$ (line 18).  
\end{itemize} 
\par To correlate the above procedure with Figure~\ref{fig:CGU}, let $p_{1}$, $p_{2}$, $p_{3}$ \& $p_{4}$ represent the pairings with highest $\mu d_j^{t}$, for $k = 5, 3, 4~ \&~ 6$, respectively. Hence, all the involved flights shall constitute $\mathcal{F}_U^t$, leading to several pairings, among which $q_1$ to $q_4$ are those with negative $\mu_j^t$s. Finally, $q_3$ and $q_4$ constitute $\mathcal{P}_{CGU}^t$, owing to their above-median $\gamma_j^t$s.
\subsection{Archiving Strategy (\texorpdfstring{$CGA$}{CGA})} \label{sec:CGA}
It may be noted that in a combinatorial optimization problem like CPOP, if a legal pairing that is part of the input for an LPP iteration, fails to be a part of the LPP solution, then it does not imply that the pairing is of no-use in absolute sense. It only implies that in that particular iteration, it did not fit among the combinations of other pairings leading to an optimal solution. Hence, it is only logical that such a pairing is re-inducted as part of the input for a subsequent iteration too. It is all the more imperative, considering that generation of legal pairing is a computationally demanding and time consuming task. Towards it, an exploratory attempt is made by the authors in \cite{aggarwal2020learning}\footnote{\cite{aggarwal2020learning} proposed a first-of-its-kind online-learning framework, based on the Variational Graph Auto-Encoders (\cite{kipf2016variational}). In that, during an intermittent LPP iteration, the learning algorithm attempts to learn plausible patterns among the flight-connections of the past LPP solutions, which are then used in the same iteration by a heuristic to generate new pairings.}. In that, though a considerable amount of cost-improvement is observed in the LPP-solutioning phases of the optimizer, the total run-time of the optimizer increased drastically. This became the motivation for the development of the \textit{archiving} strategy ($CGA$) proposed here. In any iteration $t$, $CGA$:
\begin{itemize}
	\item maintains an \textit{archive} $\mathcal{A}^t$ of the previously generated pairings, such that any pairing is stored/retrieved with reference to a unique index $(f_m,f_n)$ reserved for any legal flight-pair in that pairing. 
	\item introduces a new measure, namely, \textit{reduced cost estimator} ($\eta_{mn}^t$, Equation~\ref{eq:redCostEst}) for a flight-pair in $\mathcal{A}^t$, indexed as $(f_m,f_n)$. This estimator can be computed for all the flight-pairs present across all the pairings in the current LPP solution (by fetching $f_m$, $f_n$, $y_m$ and $y_n$). 
	\begin{flalign}
		\eta_{mn}^{t} &= flying\_cost(f_m) +  flying\_cost(f_n) - y^{t}_m - y^{t}_n = \sum_{i \in \{m,n\}} \left(flying\_cost(f_i) - y^{t}_i \right) \label{eq:redCostEst}
	\end{flalign}
	Notably, this formulation is analogous to Equation~\ref{eq:redCost}, just that instead of complete cost for a pairing, only the flying costs corresponding to the flights in a legal pair, are accounted for. Given this, $\eta_{mn}^t$ may be seen as an indicator of $\mu_j$ \textit{at the flight-pair level}. 
	\item recognizes that towards further improvement in the current LPP solution, it may be prudent to include as a part of the input for next LPP iteration -- fresh pairings comprising of flight-pairs with low $\eta_{mn}^t$. To this effect, $\mathcal{P}^t_{CGA}$ is in principle constituted by preferentially picking pairings in $\mathcal{A}^t$ comprising of flight-pairs reporting lower $\eta_{mn}^t$.
\end{itemize}
In doing so, this strategy pursues the goal of continual improvement in the objective function, while relying on the \textit{flight-pair} level information embedded in the LPP solution at any iteration, and utilizing the past computational efforts. 
%
\begin{algorithm2e}[htbp]
	\small
	\DontPrintSemicolon
	\SetKwComment{Comment}{$\triangleright$\ }{}
	\KwIn{$\mathcal{A}^{t-1},~\mathcal{P}^{t-1}_{LP} \cup \mathcal{P}^{t-1}_{CG},~Y^{t},~Th_A$}
	\KwOut{$\mathcal{P}^{t}_{CGA}$}
	\nosemic
	$\mathcal{A}^t \gets$ Update $\mathcal{A}^{t-1}$ with each pairing $p_j \in \mathcal{P}^{t-1}_{LP} \cup \mathcal{P}^{t-1}_{CG}$\; 
	$Random \gets$ Generate a random integer between 0 \& $Th_A$\;
	Compute $\eta^t_{mn}$ for each flight-pair $(f_m,f_n) \in \mathcal{A}^t$\Comment*{Use Equation~\ref{eq:redCostEst}}
	Sort elements of $\mathcal{A}^{t}$ in ascending order with-respect-to the $\eta^{t}_{mn}$s of the corresponding flight-pairs\;
	\For{each element $\left((f_m,f_n),\mathcal{P}_{mn}^t\right) \in~sorted~\mathcal{A}^t$}{
		$\mathcal{P}_A^t \gets$ Randomly select a maximum of $Random$ number of pairings from $\mathcal{P}_{mn}^t$\;
		\If{$|\mathcal{P}_A^t| \geq (Random \cdot Random)$}{
			Break\;
		}
	}
	Compute $\mu^{t}_j$ of all pairings in $\mathcal{P}_{A}^{t}$\Comment*{Use Equation~\ref{eq:redCost}}
	$\mathcal{P}^{t}_{CGA} \gets $ Select only the pairings with negative $\mu_j^t$ from $\mathcal{P}_{A}^{t}$\;
	\textbf{return} $\mathcal{P}^{t}_{CGA}$\;
	\caption{Procedure for the $CGA$ strategy, at an LPP iteration $t$}
	\label{alg:CGA}
\end{algorithm2e}
\par The implementation of this strategy at iteration $t$, has been formalized in Algorithm~\ref{alg:CGA}. In that, the first step is to obtain an updated archive, $\mathcal{A}^{t}$ (line 1), by conjoining the previous archive $\mathcal{A}^{t-1}$ and the set of pairings received as input from the previous iteration. Here, $\mathcal{A}^t$ is a dictionary, each element of which can be given by $\left((f_m,f_n),\mathcal{P}^t_{mn}\right)$, where $(f_m,f_n)$ is a legal flight-pair and $\mathcal{P}^t_{mn}$ -- the corresponding pairing set. In the first iteration ($t=1$), the archive $\mathcal{A}^{1}$ is created using the pairings from the initial feasible solution. Notably, with generation of approximately a million pairings in each LPP iteration, the size of $\mathcal{A}^t$ grows at an alarming rate. Hence, for tractability, this strategy extracts only a subset of pairings, say $\mathcal{P}_A^t$, from the current archive $\mathcal{A}^t$. The size of $\mathcal{P}_A^t$ is controlled at two levels, by $Random$ -- a random integer between $0$ and a pre-defined upper threshold, namely, $Th_A$ (the criterion for its setting is highlighted in Section~\ref{sec:parameter-setting}). In that:
\begin{itemize}
	\item firstly, $\eta_{mn}^t$ is computed for each flight-pair in $\mathcal{A}^t$ (line 3), allowing for sorting of these pairs in ascending order of their $\eta_{mn}^t$ values (line 4). Then the first $Random$ flight-pairs from the sorted list are considered (lines 2 to 5). 
	\item subsequently, for each considered flight-pair, $(f_m,f_n)$, a fixed number of pairings are selected randomly from $\mathcal{P}^t_{mn}$, and this number is also same as $Random$ (lines 6 to 9).
\end{itemize}
Once $\mathcal{P}_A^t$ is constituted as above, $\mu_j^t$ is computed (line 11) for each of its pairings, and those with negative $\mu_j^t$ values are selected to constitute the new pairing set, namely $\mathcal{P}_{CGA}^t$ (line 12).
\subsection{Random Exploration Strategy (\texorpdfstring{$CGR$}{CGR})} \label{sec:CGR}
It is critical to note that each of the CG strategy presented above: (a) \textit{deadhead reduction} ($CGR$ offering $\mathcal{P}_{CGD}^t$) pursuing optimal features at a \textit{set level}, (b) \textit{crew utilization enhancement} ($CGU$ offering $\mathcal{P}_{CGU}^t$) pursuing optimal features at a \textit{pairing level}, and (c) \textit{archiving} ($CGA$ offering $\mathcal{P}_{CGA}^t$) based on a \textit{flight-pair level} information, involves the use of  current LPP solution. In contrast, the \textit{random exploration} strategy proposed here, namely, $CGR$ aims to pursue random and unbiased exploration of the pairings' space, independent of the current LPP solution. It involves generation of new pairings for a random selected set of legal duties, of which only the pairings with negative reduced cost are selected to constitute a new pairing set, $\mathcal{P}_{CGR}^t$. Here, a random set of legal duties is used instead of a random set of flights, as the former has a higher probability of generating legal pairings, given that  satisfaction of most of the legality constraints is ensured during their pre-processing. 
\begin{algorithm2e}[htbp]
	\small
	\DontPrintSemicolon
	\SetKwComment{Comment}{$\triangleright$\ }{}
	\KwIn{$\mathcal{D}_{all},~Y^{t},~Th_R$}
	\KwOut{$\mathcal{P}^{t}_{CGR}$}
	\nosemic $Random \gets$ Generate a random integer between 0 \& $Th_R$\;
	$\mathcal{D}_{R}^t \gets$ Select $Random$ number of duties with-respect-to each crew base, from $\mathcal{D}_{all}$ without replacement\;
	Generate all possible pairings using $\mathcal{D}_{R}^t$ and compute their $\mu^{t}_{j}$\Comment*{Use Equation~\ref{eq:redCost}}
	$\mathcal{P}^{t}_{CGR} \gets$ Select only the pairings with negative $\mu^{t}_j$\;
	\textbf{return} $\mathcal{P}^{t}_{CGR}$\;
	\caption{Procedure for the $CGR$ strategy, at an LPP iteration $t$}
	\label{alg:CGR}
\end{algorithm2e}
The implementation of this strategy at iteration $t$, has been formalized in Algorithm~\ref{alg:CGR}. Its input involves the set of all legal duties $\mathcal{D}_{all}$, the dual vector $Y^{t}$, and a pre-defined threshold, denoted by $Th_R$ (the criterion for its setting is highlighted in Section~\ref{sec:parameter-setting}). A random integer between $0$ and $Th_R$, namely, $Random$, defines the number of duties with-respect-to each crew base, for which new pairings are to be generated. This cumulatively helps constitute the subset of legal duties, $\mathcal{D}^t_{R}$ (line 2), for which pairings are generated and the corresponding $\mu_j^t$ values are computed (line 3). Finally, only the pairings with negative $\mu_j^t$ values constitute $\mathcal{P}^{t}_{CGR}$ (line 4), which forms a part of the input for the next LPP iteration.\\\\
\noindent \textbf{Overarching CG heuristic based on the four CG strategies}\\
In the wake of the context set up early in this section, and the CG strategies, namely, $CGD$,~$CGU$,~$CGA$~\&~$CGR$, discussed individually, the overarching CG heuristic has been presented in the Algorithm~\ref{alg:CGheuristic}. 
\begin{algorithm2e}[htbp]
	\DontPrintSemicolon
	\small
	\KwIn{$\mathcal{D}_{all},~\mathcal{P}^{t}_{LP},~X_{LP}^{t},~Y^{t},~\mathcal{A}^{t-1},~\mathcal{P}^{t-1}_{CG} \cup \mathcal{P}^{t-1}_{LP},~Th_R,~Th_U,~Th_D,~Th_A$}
	\KwOut{$\mathcal{P}_{CG}^{t}$}
	\nosemic 
	$\mathcal{P}^{t}_{CGD} \gets CGD(\mathcal{D}_{all},\mathcal{P}^{t}_{LP},Y^{t},Th_D)$\;
	$\mathcal{P}^{t}_{CGU} \gets CGU(\mathcal{P}^{t}_{LP},X_{LP}^{t},Y^{t},Th_U)$\;
	$\mathcal{P}^{t}_{CGA} \gets CGA(\mathcal{A}^{t-1},\mathcal{P}^{t-1}_{CG} \cup \mathcal{P}^{t-1}_{LP},Y^{t},Th_A)$\;
	$\mathcal{P}^{t}_{CGR} \gets CGR(\mathcal{D}_{all},Y^{t},Th_R)$\;
	$\mathcal{P}_{CG}^{t} \gets \mathcal{P}^{t}_{CGD} \cup \mathcal{P}^{t}_{CGU} \cup \mathcal{P}^{t}_{CGA} \cup \mathcal{P}^{t}_{CGR}$\;
	\textbf{return} $\mathcal{P}_{CG}^{t}$\;
	\caption{Procedure for the proposed CG heuristic, at an LPP iteration $t$}
	\label{alg:CGheuristic}
\end{algorithm2e}
Notably, the input to the procedure comprises of the set of all legal duties $\mathcal{D}_{all}$, the current LPP solution $\mathcal{P}^{t}_{LP}$, its corresponding primal vector $X_{LP}^{t}$ \& dual vector $Y^{t}$, the previous archive $\mathcal{A}^{t-1}$, the pairing set input from the previous iteration $\mathcal{P}^{t-1}_{LP} \cup \mathcal{P}^{t-1}_{CG}$, and the pre-defined thresholds, namely, $Th_D$, $Th_U$, $Th_A$, \& $Th_R$, for their corresponding CG strategies. It is also evident (lines 1 to 5) that at any given LPP iteration, these CG strategies are executed in a sequential manner, each leading up to a part of the input for the next LPP iteration ($\mathcal{P}^{t}_{CG}$).
\section{Computational Experiments} \label{sec:exp}
This section presents the test cases used to investigate the performance of the proposed CG heuristic and its underlying strategies. The computational setup and the parameter settings used for this investigation are also highlighted followed by discussion of the experimental results.
\subsection{Test Cases} \label{sec:exp-testcases}
The real-world airline test cases used for experimentation are cited in Table~\ref{tab:testcases}. Each of these test cases involves a weekly flight-schedule, and have been provided by the research consortium's Industrial Sponsor. 
\begin{table}[htbp]
	\small 
	\caption{Real-world airline test cases used in this research work}
	\begin{center}
		\begin{tabular}{cccc}
			\toprule
			\textbf{Test Cases} & \textbf{$\#$Flights} & \textbf{$\#$Crew Bases} & \textbf{$\#$Legal Duties} \\ 				\midrule
			TC-1 & 3202 & 15 & 454205 \\
			TC-2 & 3228 & 15 & 464092 \\
			TC-3 & 3229 & 15 & 506272 \\
			TC-4 & 3265 & 15 & 446937 \\
			TC-5 & 4212 & 15 & 737184 \\
			\bottomrule
		\end{tabular}
		\label{tab:testcases}
	\end{center}
\end{table}
\begin{figure}[htbp]
	\centering
	\begin{subfigure}[t]{0.225\textwidth}
		\centering{\includegraphics[width=\columnwidth, keepaspectratio]{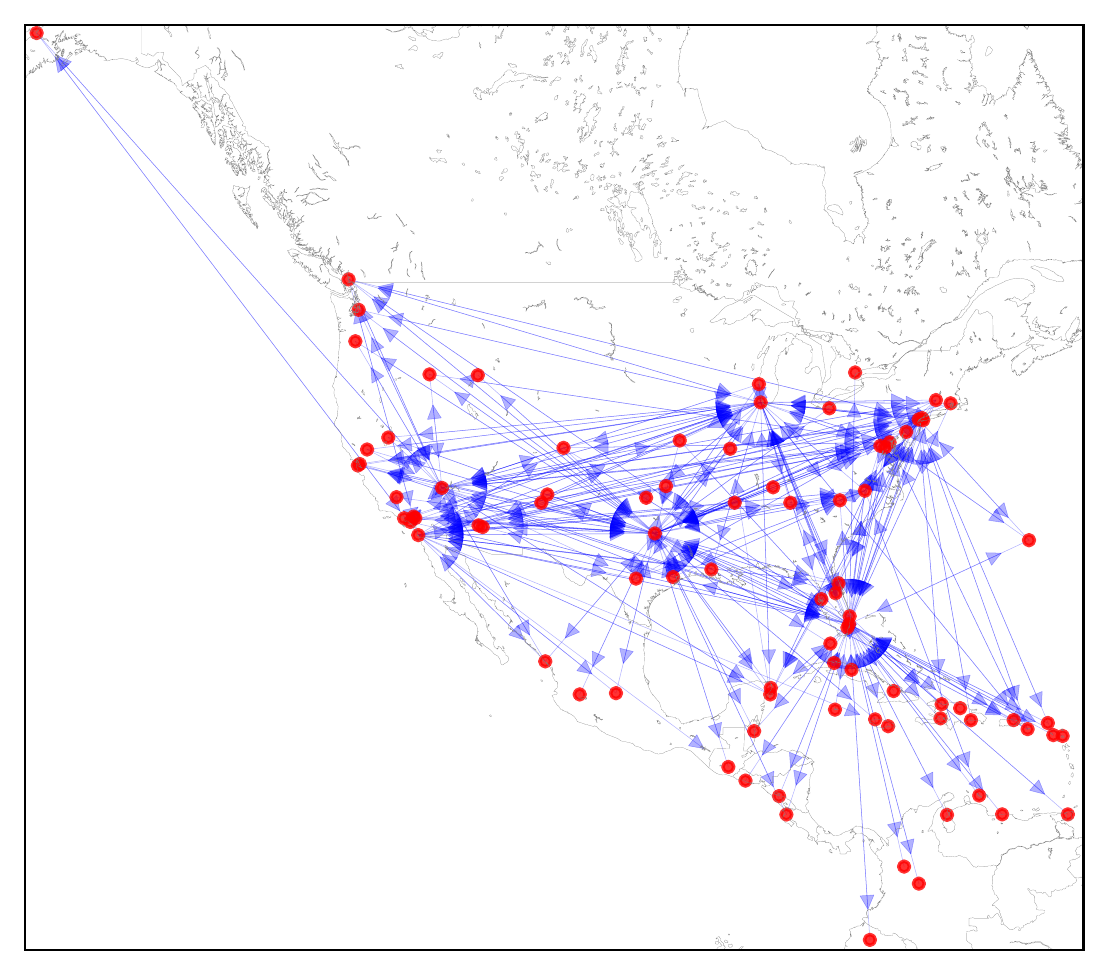}}
		\centering \caption{}
		\label{fig:fullFlightMap}
	\end{subfigure}
	\begin{subfigure}[t]{0.2475\textwidth}
		\centering{\includegraphics[width=\columnwidth, keepaspectratio]{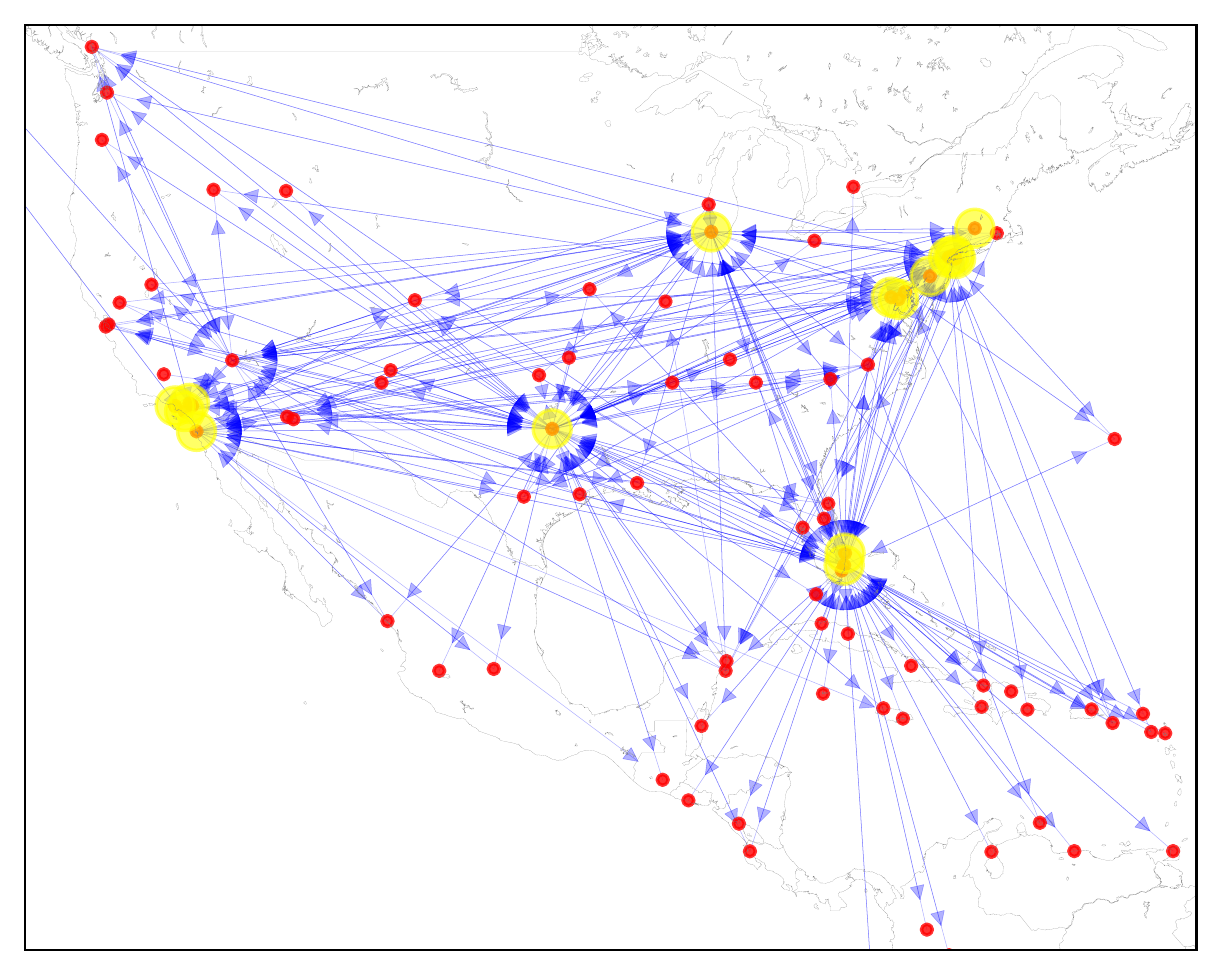}}
		\centering \caption{}
		\label{fig:flightMapWithCB}
	\end{subfigure}
	\begin{subfigure}[t]{0.2475\textwidth}
		\centering{\includegraphics[width=\columnwidth, keepaspectratio]{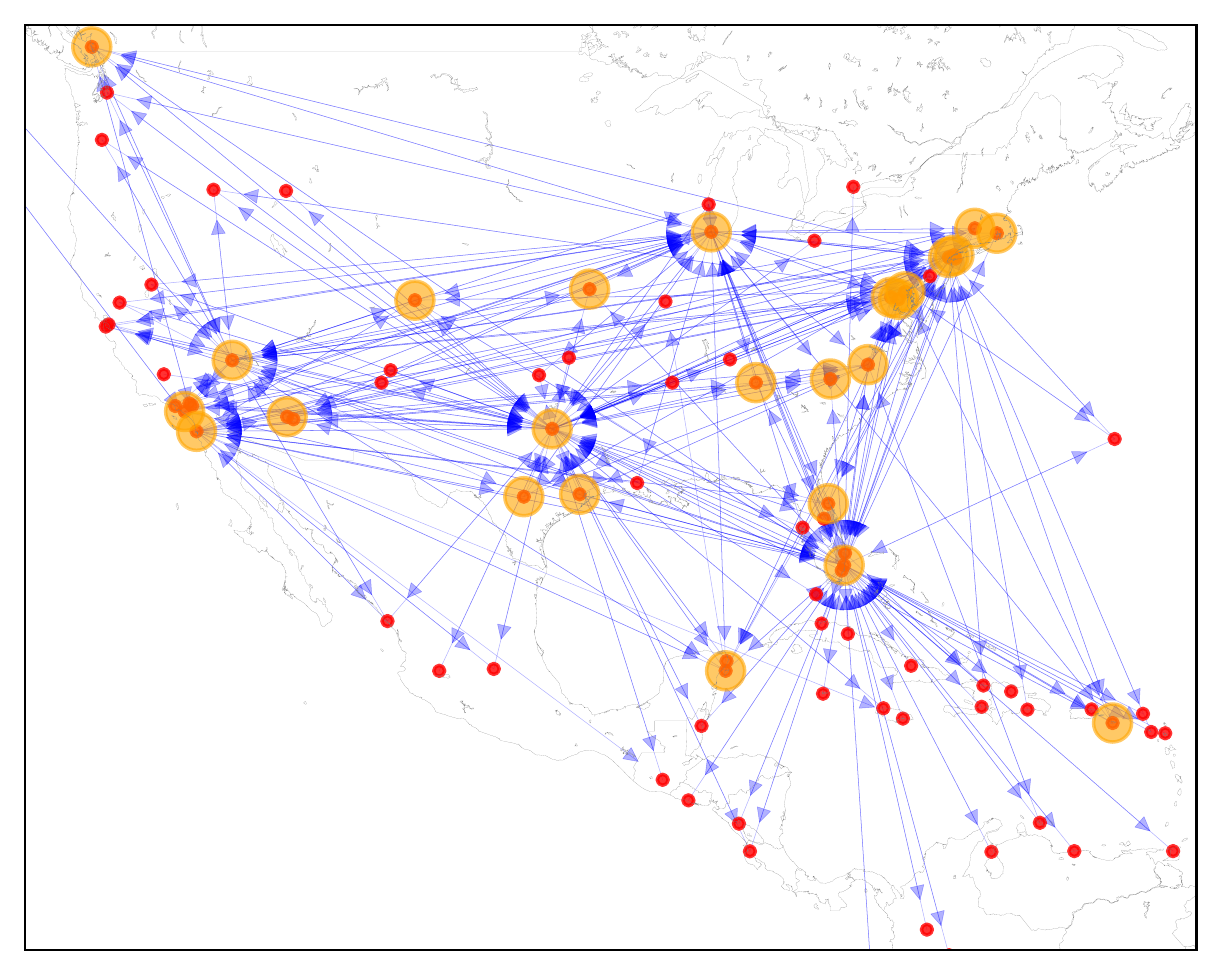}}
		\centering \caption{}
		\label{fig:flightMapWithHubs}
	\end{subfigure}
	\begin{subfigure}[t]{0.22
	\textwidth}
		\centering{\includegraphics[width=\columnwidth, keepaspectratio]{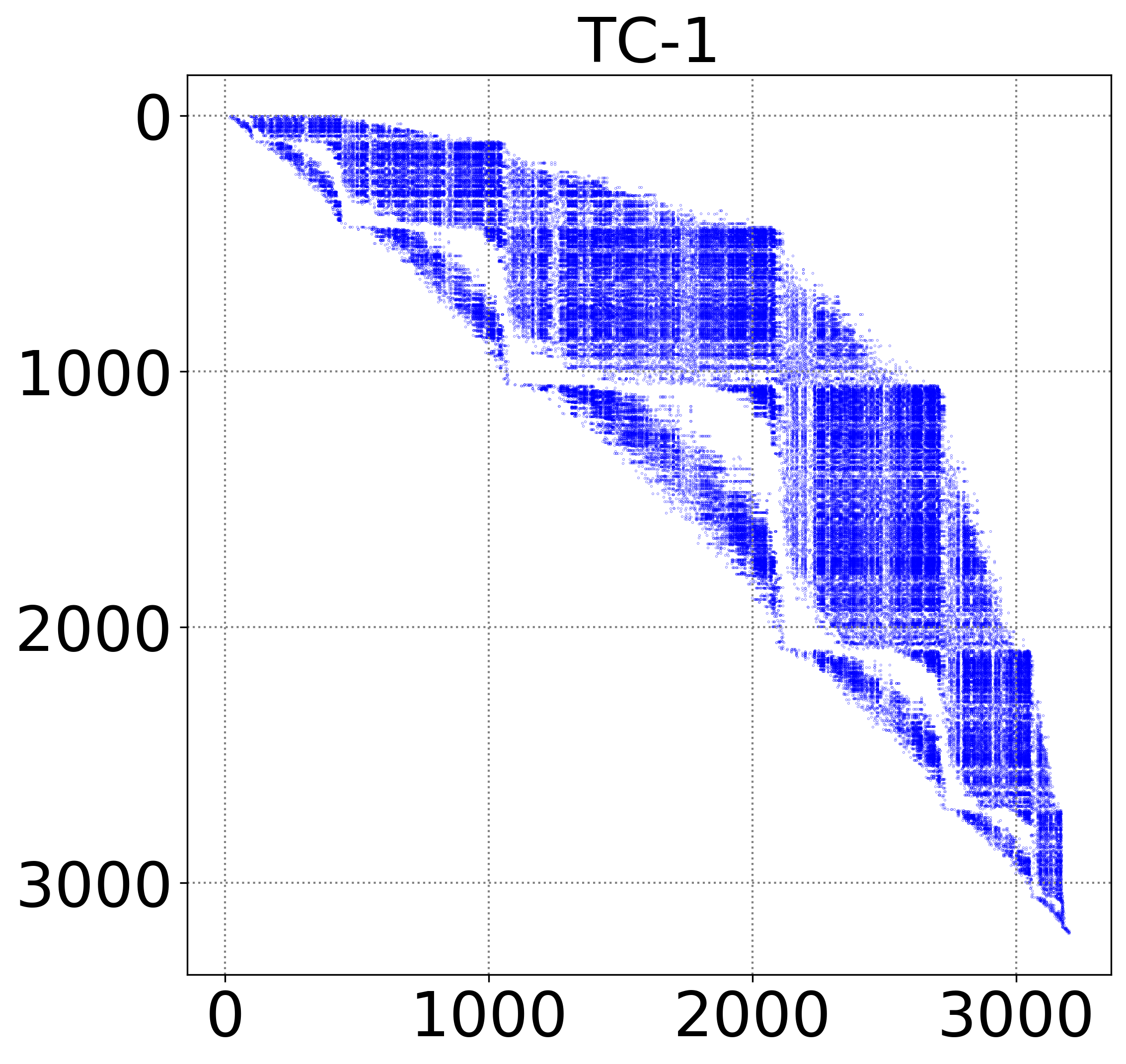}}
		\centering \caption{}
		\label{fig:plotTCs}
	\end{subfigure}
	\caption{(a) Geographical representation of TC-1's flight network, where the red nodes and blue edges represent the airports and scheduled flights, respectively; (b) Zoomed view of TC-1's flight network with all 15 crew bases highlighted in yellow color; (c) Zoomed view of TC-1's flight network with major hub airports (serving $\geq$4 other airports) highlighted in orange color; and (d) \textit{legal} flight-connections, each represented by a point in the plot, where for a flight marked on the y-axis, the connecting flight is marked on the x-axis. Notably, in (d), the flights involved in TC-1 are indexed in the ascending order of their departure time-stamps.}
	\label{fig:TCs}
\end{figure}
The columns in Table~\ref{tab:testcases}, in order of their occurrence, highlight the notations for the different test cases; the number of its constituent flights; the number of constituent crew bases; and the total number of legal duties involved, respectively. It is critical to recognize that the \textit{challenge associated with the solutioning of these test cases is attributed to the fact that an inestimably large number of legal pairings are possible for them. This explosion of possible legal pairings is due to not just the large number of flights involved, but also the complexity of the flight networks involved (characterized by multiple crew bases and multiple hub-and-spoke sub-networks)}, in that, the number of possible legal pairings:
\begin{itemize}
    \item may grow linearly with the number of crew bases, as they provide the starting and end point for a possible pairing.
    \item may grow exponentially with the number of hubs involved, as each could serve as an intermediary airport, allowing more flight connections and coverage of more airports within a given duty \textit{elapsed time}, as opposed to a network with a single hub. This is pertinent particularly when the objective of CPOP encompasses minimization of \textit{excess pay} which in turn links with maximal utilization of duty \textit{elapsed time}.
\end{itemize}
\par As a sample instance, the geographical representation of the flight network associated with TC-1, its complexity (the presence of multiple crew bases and multiple hub-and-spoke sub-networks) and the \textit{legal} flight connections involved, are illustrated in Figure~\ref{fig:TCs}. In that, the full-scale view of the underlying flight network is portrayed by Figure~\ref{fig:fullFlightMap}. Moreover, Figures~\ref{fig:flightMapWithCB}~and~\ref{fig:flightMapWithHubs} illustrates a zoomed view of the flight network with all 15 crew bases (highlighted in yellow color) and major hub airports (highlighted in orange color), respectively. 
Furthermore, the pattern visible in Figure~\ref{fig:plotTCs} could be attributed to the (minimum and maximum) limits on the sit-time and overnight-rest constraints. For instance, a flight, say $f_{500}$, has legal connections only with those flights that depart from the arrival airport of $f_{500}$, and whose departure-time gap (difference between its departure-time and the arrival time of $f_{500}$) lies within the minimum \& maximum allowable limit, of the sit-time or the overnight-rest.
\subsection{Computational Setup} \label{sec:computational-setup}
All the experiments are performed on an HP Z640 Workstation, powered by two Intel$^\circledR$ Xeon$^\circledR$ E5-2630v3 processors, each with 16 cores at 2.40GHz, and 96GB RAM. The codes for \textit{AirCROP} built around the proposed CG heuristic, are developed using the Python scripting language (aligned with the Industrial sponsor's larger vision and preference). In that, \textit{AirCROP} serves as a testbed to investigate the performance of the proposed CG heuristic, and all its features \& associated parameters are uniformly maintained for the experimentation. Furthermore: 
\begin{itemize}
	\item In the LPP-solutioning module, the interior-point method from Gurobi Optimizer 8.1.1 \citep{gurobi} is used to solve the primal of the LPP formulation, and the interior-point method \citep{andersen2000mosek} from SciPy library \citep{scipy} is used to solve the subsequent dual of the LPP formulation. Moreover, the termination of the LPP-solutioning is enforced, if the cost improvement between two successive LPP iterations falls below 100 USD and this trend is successively repeated 10 times. These limits are set to achieve an LPP solution with a sufficiently good cost in a reasonably good time.
	\item The MIP solver of Gurobi Optimizer 8.1.1 is used in the IPP-solutioning. As mentioned before, the branch-and-bound search on a large-scale IPP is time-intensive. Hence, the IPP-solutioning module is allowed to run for a time limit of 20 minutes, if not terminated by itself.
\end{itemize}   
In principle, an \textit{AirCROP}-run terminates when the cost of the IPP solution matches the cost of its input LPP solution. However, for practical considerations on the time limit, if the IPP and LPP costs do not conform with each other despite 30 LPP-IPP interactions, or up to 30 hours of total run time, the termination of \textit{AirCROP} is enforced. 
\subsection{Parameter Settings} \label{sec:parameter-setting}
The thresholds $Th_D,~Th_U,~Th_A~\&~Th_R$, used in their respective CG strategies, regulate the size of $\mathcal{P}_{CG}^{t}$. To identify their best settings, it is important to understand the effect of the size of $\mathcal{P}^t_{CG}$ on the LPP solution's cost quality \& required run-time, and the demand on the computational resources (dominantly, RAM), as highlighted below: 
\begin{itemize}
	\item for a relatively small-sized $\mathcal{P}^t_{CG}$, the alternatives available in pairings to ensure further cost improvement shall be restricted, leading to smaller cost benefits in each phase of LPP-solutioning. This would necessitate far more LPP-IPP interactions, to reach the near-optimal cost. This by itself is not a challenge, however, significant amount of additonal run time may be required , since: (a) each phase of LPP-solutioning implies a minimum of 10 LPP iterations, before it could be terminated, (b) such unpromising phases when invoked repeatedly, may consume significant run-time, yet, without reasonable cost benefit. 
	\item On the other hand, for a very large-sized $\mathcal{P}_{CG}^{t}$, though the potential for significant cost benefits may exist, the demand on the RAM may become overwhelming for any LPP solutioning phase to proceed. 
\end{itemize}
The above considerations suggest that the size of $\mathcal{P}_{CG}^{t}$ may neither be too small nor too large. Factoring these, the experiments here aim at $\mathcal{P}_{CG}^{t}$ sized approximately of a million pairings (significant size, yet, not overwhelming for 96 GB RAM). Furthermore, for a search that is not biased in favour of any particular CG strategy constituting the heuristic, it is also important that the number of pairings from each strategy are more or less equable. To ensure this, the random integers generated in each CG strategy are selected from the same half of their respective ranges. For instance, if a random integer in first CG strategy is generated from the upper/lower half of its respective range, then so are the random integers in the remaining three CG strategies. 
\subsection{Results \& Observations} \label{results}
This section, primarily has three goals. The first goal is to be investigate if the use of the proposed CG heuristic within \textit{AirCROP} offers any advantage over the stand-alone use of just the $CGR$ strategy, or, just the standard CG technique (with exact-pricing strategy, as mentioned in Section~\ref{sec:relatedwork}) alone. In doing so, the endeavor is to assess whether the quest of the proposed CG heuristic to balance the \textit{random search} and \textit{guided search} (via exploitation of optimal solution features) infuses better search efficiency compared to the \textit{random-search} alone, or, the \textit{exact-search} alone. The second goal is to investigate the impact of the variability (varying pseudo-random numbers' seed) on the performance of the proposed CG heuristic vis-$\grave{a}$-vis the stand-alone use of just the $CGR$ strategy, or, just the standard CG technique. The final goal is to investigate if the proposed CG heuristic based on four CG strategies, always performs better compared to any other possible combination comprising of only two or three strategies, at a time. In doing so, the endeavor is to assess how critical the combined interplay of the four CG strategies is, compared to that offered by any subset (sized two or three) of these strategies.
\subsubsection{\textit{AirCROP} performance: Proposed CG heuristic vis-\texorpdfstring{$\grave{a}$}{a'}-vis CG Random and CG Standard}
While pursuing the first goal, the proposed CG heuristic is compared with two other CG mechanisms: (a) the stand-alone use of just the $CGR$ strategy, and (b) the standard CG technique with exact-pricing strategy. For a fair comparison of the proposed CG heuristic with the former CG mechanism, it is imperative that the number of pairings considered while employing $CGR$ alone are equable to those considered by the proposed CG heuristic. Towards it, $Th_R$ value used for stand-alone-$CGR$ is so increased that the resulting number of pairings approximately size up to those resulting from the proposed CG heuristic (which still includes $Th_R$ based $CGR$ as a constituent strategy). For the results of the above comparison (presented in Table~\ref{tab:allCGs}), the stand-alone use of $CGR$ strategy with an increased threshold is referred to as \textit{CG Random}, and the stand-alone use of standard CG technique is referred to as \textit{CG Standard}. While, the standard CG technique is discussed in Section~\ref{sec:relatedwork}, its details are shared in the subsequent text.
\afterpage{%
	\clearpage
	\thispagestyle{empty}
	\begin{landscape}
		\centering 
		\begin{table}[htbp]
			\scriptsize
			\centering
			\caption{\textit{AirCROP} performance$^\ast$: Proposed CG heuristic vis-$\grave{a}$-vis CG Random and CG Standard}
			\resizebox{\columnwidth}{!}{%
				\begin{tabular}{M{5mm}M{10mm}M{5mm}M{10.5mm}M{7.5mm}M{7.5mm}M{7.5mm}M{7.5mm}M{7.5mm}M{7.5mm}M{7.5mm}M{7.5mm}M{7.5mm}M{7.5mm}M{7.5mm}M{7.5mm}M{7mm}M{7mm}M{7mm}M{7mm}M{7mm}M{7mm}M{9mm}}
					\toprule
					\multirow{3}{*}{\textbf{Test}} & \multicolumn{2}{c}{\multirow{3}{*}{\textbf{CG}}} & \multirow{4}{*}{\bm{$\mathcal{P}_{IFS}$}} & \multicolumn{18}{c}{\textbf{LPP-IPP Interactions (\bm{$T, \mathcal{P}_{LP}^T / \mathcal{P}_{IP}^T$})}} & \multirow{3}{*}{\textbf{Final}} \\
					\cmidrule(lr){5-22}
					\multirow{3}{*}{\textbf{Case}} & \multicolumn{2}{c}{\multirow{3}{*}{\textbf{mechanism}}} &  & \multicolumn{2}{c}{\textbf{1}} & \multicolumn{2}{c}{\textbf{2}} & \multicolumn{2}{c}{\textbf{3}} & \multicolumn{2}{c}{\textbf{4}} & \multicolumn{2}{c}{\textbf{5}} & \multicolumn{2}{c}{\textbf{6}} & \multicolumn{2}{c}{\textbf{7}} & \multicolumn{2}{c}{\textbf{8}} & \multicolumn{2}{c}{\textbf{9}} & \multirow{3}{*}{\textbf{Result}} \\
					\cmidrule(lr){5-6} \cmidrule(lr){7-8} \cmidrule(lr){9-10} \cmidrule(lr){11-12} \cmidrule(lr){13-14} \cmidrule(lr){15-16} \cmidrule(lr){17-18} \cmidrule(lr){19-20} \cmidrule(lr){21-22}
					& & &  & \bm{$\mathcal{P}_{LP}^1$} & \bm{$\mathcal{P}_{IP}^1$} & \bm{$\mathcal{P}_{LP}^2$} & \bm{$\mathcal{P}_{IP}^2$} & \bm{$\mathcal{P}_{LP}^3$} & \bm{$\mathcal{P}_{IP}^3$} & \bm{$\mathcal{P}_{LP}^4$} & \bm{$\mathcal{P}_{IP}^4$} & \bm{$\mathcal{P}_{LP}^5$} & \bm{$\mathcal{P}_{IP}^5$} & \bm{$\mathcal{P}_{LP}^6$} & \bm{$\mathcal{P}_{IP}^6$} & \bm{$\mathcal{P}_{LP}^7$} & \bm{$\mathcal{P}_{IP}^7$} & \bm{$\mathcal{P}_{LP}^8$} & \bm{$\mathcal{P}_{IP}^8$} & \bm{$\mathcal{P}_{LP}^9$} & \bm{$\mathcal{P}_{IP}^9$} &  \\
					\midrule
					\multirow{9}{*}{\textbf{TC-1}} & \textbf{CG} & \textbf{Cost} & 74945982 & 3465379 & 3689312 & 3467054 & 3587416 & 3467182 & 3547549 & 3467289 & 3518776 & 3470404 & 3481129 & 3470142 & 3478695 & 3471975 & 3471975 &  &  &  &  & \textbf{3471975} \\
					\cmidrule(lr){3-23}
					& \textbf{heuristic} & \textbf{Time} & 00:04 & 04:07 & 00:20 & 01:31 & 00:20 & 01:29 & 00:20 & 01:19 & 00:05 & 00:45 & 00:01 & 00:43 & 00:01 & 00:38 & 00:01 &  &  &  &  & \textbf{11:44} \\
					\cmidrule(lr){2-23}
					& \textbf{CG} & \textbf{Cost} & 74945982 & 3475872 & 3858042 & 3477979 & 3741546 & 3478246 & 3646239 & 3484971 & 3558809 & 3534145 & 3534145 &  &  &  &  &  &  &  &  & \textbf{3534145} \\
					\cmidrule(lr){3-23}
					& \textbf{Random} & \textbf{Time} & 00:04 & 04:34 & 00:20 & 03:52 & 00:20 & 03:54 & 00:20 & 02:20 & 00:20 & 00:39 & 00:01 &  &  &  &  &  &  &  &  & \textbf{16:44} \\
					\cmidrule(lr){2-23}
					& \textbf{CG} & \textbf{Cost} & 74945982 & 3581668 & 4036848 & 3946632 & 3946632 &  &  &  &  &  &  &  &  &  &  &  &  &  &  & \textbf{3946632} \\
					\cmidrule(lr){3-23}
					& \textbf{Standard} & \textbf{Time} & 00:04 & 06:44 & 00:20 & 00:28 & 00:01 &  &  &  &  &  &  &  &  &  &  &  &  &  &  & \textbf{07:37} \\
					\midrule
					\multirow{9}{*}{\textbf{TC-2}} & \textbf{CG} & \textbf{Cost} & 129221508 & 3495054 & 3779111 & 3492575 & 3632511 & 3494078 & 3573900 & 3494145 & 3532658 & 3495669 & 3529639 & 3496165 & 3514831 & 3497105 & 3501399 & 3496510 & 3513764 & 3497010 & 3497010 & \textbf{3497010} \\
					\cmidrule(lr){3-23}
					& \textbf{heuristic} & \textbf{Time} & 00:08 & 03:58 & 00:20 & 03:34 & 00:20 & 01:09 & 00:20 & 01:22 & 00:01 & 00:55 & 00:01 & 01:10 & 00:01 & 00:52 & 00:01 & 00:40 & 00:01 & 00:37 & 00:01 & \textbf{15:31} \\
					\cmidrule(lr){2-23}
					& \textbf{CG} & \textbf{Cost} & 129221508 & 3511158 & 3958787 & 3506421 & 3747013 & 3510784 & 3697667 & 3508414 & 3635284 & 3515155 & 3597757 & 3517489 & 3558292 & 3536578 & 3536578 &  &  &  &  & \textbf{3536578} \\
					\cmidrule(lr){3-23}
					& \textbf{Random} & \textbf{Time} & 00:08 & 02:40 & 00:20 & 03:47 & 00:20 & 02:08 & 00:20 & 04:03 & 00:20 & 01:30 & 00:20 & 01:58 & 00:01 & 00:36 & 00:01 &  &  &  &  & \textbf{18:32} \\
					\cmidrule(lr){2-23}
					& \textbf{CG} & \textbf{Cost} & 129221508 & 3600319 & 4076330 & 3909212 & 3911839 & 3902100 & 3902100 &  &  &  &  &  &  &  &  &  &  &  &  & \textbf{3902100} \\
					\cmidrule(lr){3-23}
					& \textbf{Standard} & \textbf{Time} & 00:08 & 10:03 & 00:20 & 00:32 & 00:01 & 00:23 & 00:01 &  &  &  &  &  &  &  &  &  &  &  &  & \textbf{11:28} \\
					\midrule
					\multirow{9}{*}{\textbf{TC-3}} & \textbf{CG} & \textbf{Cost} & 40721754 & 3485462 & 3730120 & 3486563 & 3615698 & 3486182 & 3559567 & 3486190 & 3528763 & 3495213 & 3503498 & 3495142 & 3495142 &  &  &  &  &  &  & \textbf{3495142} \\
					\cmidrule(lr){3-23}
					& \textbf{heuristic} & \textbf{Time} & 00:03 & 05:08 & 00:20 & 02:36 & 00:20 & 02:09 & 00:20 & 02:22 & 00:02 & 01:31 & 00:01 & 00:47 & 00:01 &  &  &  &  &  &  & \textbf{15:40} \\
					\cmidrule(lr){2-23}
					& \textbf{CG} & \textbf{Cost} & 40721754 & 3495505 & 3879551 & 3495881 & 3716709 & 3496945 & 3616849 & 3499394 & 3568934 & 3499705 & 3548203 & 3524914 & 3524914 &  &  &  &  &  &  & \textbf{3524914} \\
					\cmidrule(lr){3-23}
					& \textbf{Random} & \textbf{Time} & 00:03 & 05:04 & 00:20 & 05:18 & 00:20 & 04:20 & 00:20 & 04:11 & 00:20 & 07:22 & 00:08 & 00:45 & 00:01 &  &  &  &  &  &  & \textbf{28:32} \\
					\cmidrule(lr){2-23}
					& \textbf{CG} & \textbf{Cost} & 40721754 & 3616678 & 4117779 & 3991375 & 3991375 &  &  &  &  &  &  &  &  &  &  &  &  &  &  & \textbf{3991375} \\
					\cmidrule(lr){3-23}
					& \textbf{Standard} & \textbf{Time} & 00:03 & 09:47 & 00:20 & 00:40 & 00:01 &  &  &  &  &  &  &  &  &  &  &  &  &  &  & \textbf{10:51} \\
					\midrule
					\multirow{9}{*}{\textbf{TC-4}} & \textbf{CG} & \textbf{Cost} & 93198213 & 3593922 & 3812726 & 3595165 & 3698944 & 3597266 & 3649587 & 3597369 & 3624450 & 3597846 & 3616136 & 3598890 & 3605885 & 3599384 & 3601844 & 3599539 & 3599539 &  &  & \textbf{3599539} \\
					\cmidrule(lr){3-23}
					& \textbf{heuristic} & \textbf{Time} & 00:05 & 04:01 & 00:20 & 02:07 & 00:20 & 01:21 & 00:01 & 01:21 & 00:01 & 00:54 & 00:01 & 01:09 & 00:01 & 00:36 & 00:01 & 00:30 & 00:01 &  &  & \textbf{12:50} \\
					\cmidrule(lr){2-23}
					& \textbf{CG} & \textbf{Cost} & 93198213 & 3603090 & 3864971 & 3606451 & 3809623 & 3606551 & 3749713 & 3609458 & 3677935 & 3651825 & 3653528 & 3643191 & 3643191 &  &  &  &  &  &  & \textbf{3643191} \\
					\cmidrule(lr){3-23}
					& \textbf{Random} & \textbf{Time} & 00:05 & 05:16 & 00:20 & 02:19 & 00:20 & 02:45 & 00:20 & 02:28 & 00:20 & 00:59 & 00:01 & 01:14 & 00:01 &  &  &  &  &  &  & \textbf{16:28} \\
					\cmidrule(lr){2-23}
					& \textbf{CG} & \textbf{Cost} & 93198213 & 3707603 & 4217994 & 4093020 & 4093020 &  &  &  &  &  &  &  &  &  &  &  &  &  &  & \textbf{4093020} \\
					\cmidrule(lr){3-23}
					& \textbf{Standard} & \textbf{Time} & 00:05 & 04:58 & 00:20 & 00:24 & 00:01 &  &  &  &  &  &  &  &  &  &  &  &  &  &  & \textbf{05:48} \\
					\midrule
					\multirow{9}{*}{\textbf{TC-5}} & \textbf{CG} & \textbf{Cost} & 131443284 & 4584526 & 4997918 & 4590082 & 4795955 & 4589793 & 4743617 & 4592119 & 4685270 & 4591363 & 4672393 & 4591374 & 4631706 & 4593502 & 4604572 & 4593256 & 4593256 &  &  & \textbf{4593256} \\
					\cmidrule(lr){3-23}
					& \textbf{heuristic} & \textbf{Time} & 00:09 & 07:40 & 00:20 & 03:16 & 00:20 & 03:45 & 00:20 & 02:41 & 00:20 & 02:19 & 00:20 & 02:26 & 00:01 & 01:28 & 00:01 & 01:06 & 00:01 &  &  & \textbf{26:33} \\
					\cmidrule(lr){2-23}
					& \textbf{CG} & \textbf{Cost} & 131443284 & 4599464 & 5140530 & 4617224 & 4972440 & 4629868 & 4800427 & 4749477 & 4749477 &  &  &  &  &  &  &  &  &  &  & \textbf{4749477} \\
					\cmidrule(lr){3-23}
					& \textbf{Random} & \textbf{Time} & 00:09 & 10:09 & 00:20 & 04:42 & 00:20 & 06:10 & 00:20 & 01:02 & 00:01 &  &  &  &  &  &  &  &  &  &  & \textbf{23:13} \\
					\cmidrule(lr){2-23}
					& \textbf{CG} & \textbf{Cost} & 131443284 & 4802101 & 5558634 & 4814379 & 5279814 & 4836901 & 5109643 &  &  &  &  &  &  &  &  &  &  &  &  & \textbf{5109643} \\
					\cmidrule(lr){3-23}
					& \textbf{Standard} & \textbf{Time} & 00:09 & 10:06 & 00:20 & 10:14 & 00:20 & 10:11 & 00:20 &  &  &  &  &  &  &  &  &  &  &  &  & \textbf{31:40}\\
					\bottomrule
			\end{tabular}}
			\label{tab:allCGs}
			\\\footnotesize{$^\ast$All values in the ``Cost'' rows are in USD, and all the corresponding real values are rounded-off to the next integer values. All values in the ``Time'' rows are in HH:MM, and all the corresponding seconds' values are rounded-off to the next minute values.}
		\end{table}
	\end{landscape}
	\clearpage
}
\par The latest instance that best describes the application of standard CG technique to large-scale airline CPOPs is \citet{desaulniers2020dynamic}, and the same approach has been implemented in the CG Standard here. In that, there exists a pricing sub-problem for each combination of a day of the planning horizon and a crew base, allowing the generation of only those pairings which start on the corresponding day from the corresponding crew base and end at the same crew base. Each of these pricing sub-problems is modeled as SPPRC, and is defined on an acyclic time-space network, wherein each node represents a combination of an airport and a time, and each arc represents a movement in space and/or time (for e.g. a flight, a duty, a sit-time connection, etc.). To solve it, the authors used a dynamic programming and labeling algorithm, proposed by \citet{irnich2005shortest}. For large-scale but less complex flight networks, \citet{desaulniers2020dynamic} adopted a flight-based time-space network, originally presented by \citet{saddoune2013aircrew}, wherein the arcs represent only those flights that can be covered by a pairing starting on the day associated with the respective pricing sub-problem. In this type of structure, all source-to-sink paths represent a pairing that starts and ends at the same crew base, and respects some pairing-legality constraints such as maximum duration of a pairing, etc. The other legality constraints are modeled using \textit{resource constraints} \citep{irnich2005shortest}, wherein a resource is a quantity that varies along a path and is allowed to take values only within a defined resource window at each visited node. Some examples of these resources are the number of duties in a pairing, duration of a duty, number of flights allowed in a duty, etc. However, for much larger and complex flight networks (complexity being similar to the ones tackled in this research work), \citet{desaulniers2020dynamic} used a duty-based time-space network to define each pricing sub-problem, and solved it using a heuristic labeling algorithm with a label dominance rule \citep{irnich2005shortest}. This heuristic labeling algorithm helps in eliminating the need for large number of resources to account for large number of pairing legality constraints. Here, the same approach has been adopted to implement CG Standard. In this approach, a complete set of legal duties is enumerated a priori, which is then used to define the acyclic time-space network for each pricing sub-problem. Within this network structure, each node represents a combination of an airport and a time, and each arc represents a duty from the enumerated set that can be covered by a pairing starting on the day associated with the respective pricing sub-problem.
\par As evident, Table~\ref{tab:allCGs} presents the results for all the five test cases (TC-1 to TC-5). In that, for each test case:
\begin{itemize}
	\item the column marked by ``$\mathcal{P}_{IFS}$'' highlights the cost associated with the $IFS$ that triggers the respective \textit{AirCROP}-run and the time consumed in its generation.
	\item the subsequent columns present the results of final solutions of the LPP-IPP interactions, marked by $T$. In that, the column marked by ``$\mathcal{P}_{LP}^T$'' highlights the cost associated with the final solution of the respective LPP-solutioning phase and its run-time. Similarly, the column marked by ``$\mathcal{P}_{IP}^T$'' highlights the cost associated with the final solution of the respective IPP-solutioning phase and its run-time.
	\item the last (emboldened) column marked by ``Final Solution'' highlights the results of the final CPOP/IPP solutions of the respective \textit{AirCROP}-runs, in terms of cost and total run-time consumed.
	%
\end{itemize}
Notably, from the tabulated results, the following observations can be drawn for all test cases: 
\begin{enumerate}
	\item The runs with the proposed CG heuristic led to the final IPP solutions with lower cost than those of the runs with CG Random. These cost improvements vary from 0.84\% (29,772 USD for TC-3) to 3.29\% (156,221 USD for TC-5), which if saved could translate to annual savings worth millions of dollars. A similar trend is observed between the total run-times of the two runs for all test cases, except for TC-5. These run-time improvements vary from 16.28\% (03:01 for TC-2) to 45.09\% (12:52 for TC-3).
	\item The runs with the proposed CG heuristic led to the final IPP solutions with significantly lower cost ($>$10\%) than those of the runs with CG Standard. These cost improvements vary from 10.11\% (516,387 USD for TC-5) to 12.43\% (496,233 USD for TC-3), promising significant annual savings. 
	In terms of the run-time, CG Standard bears two trends. In the case of relatively smaller test cases (TC-1 to TC-4 with $<$3300 flights), it seems to lead to premature termination, while for the larger test case (TC-5 with $>$4200 flights) it can be seen to terminate after spanning the upper limit of 30 hours. These trends could be attributed to the size of the IPP solution at $T=1$ (i.e., $|\mathcal{P}_{IP}^{1}|$), since it serves as an IFS for the LPP-solutioning phase at $T=2$. This is endorsed by the fact that for TC-1 to TC-4, $|\mathcal{P}_{IP}^{1}| \approx$ 600, while for TC-5, $|\mathcal{P}_{IP}^{1}| \approx $ 1000. Given the correlation between the size of $\mathcal{P}_{IP}^{1}$ and the timing of termination, it is plausible to argue that the premature termination in the case of TC-1 to TC-4 coincides with the inability of CG Standard to generate promising pairings with respect to the limited set in $\mathcal{P}_{IP}^{1}$.
	\item For all test cases, experiments with CG Standard reveal that the costs obtained after completion of LPP-solutioning phases, keep getting worse over the phases. In other words, $\mathcal{P}^T_{LP}$ keeps getting worse with increasing $T$. However, this trend is broken in the case of the proposed CG heuristic and CG Random. For instance: (i) in the case of the proposed CG heuristic, the cost of $\mathcal{P}^6_{LP}$ gets better than the cost of $\mathcal{P}^5_{LP}$ for TC-1, and (ii) in the case of CG Random, the cost of $\mathcal{P}^2_{LP}$ gets better than the cost of $\mathcal{P}^1_{LP}$ for TC-2.
\end{enumerate}
The above observations indicate that the search efficiency of the proposed CG heuristic is better than that of CG Random as well as CG Standard. The plausible reason is that the proposed CG heuristic not only explores the new pairings' space in a random and unbiased manner using $CGR$ (which CG Random does), but also exploits the optimal solution's features at a set-level using $CGD$, and at a pairing-level using $CGU$, while also exploiting the information on promising flight-pairs using $CGA$.
\subsubsection{Variability in \textit{AirCROP} Performance: Proposed CG heuristic vis-\texorpdfstring{$\grave{a}$}{a'}-vis CG Random and CG Standard}
In this section, empirical results to investigate the impact of variability (varying pseudo-random numbers' seed) on the performance of the proposed CG heuristic vis-$\grave{a}$-vis CG Random and CG Standard, are presented.
\par \citet{koch2011miplib} demonstrated that when the mathematical programming-based solution approaches are employed, such as the ones in \textit{AirCROP}, some variability in their final performance is observed, which is rather inevitable. \citet{lodi2013performance} noted that the factors such as changing the floating-point arithmetic, permuting the constraints/variables of the underlying mathematical models, changing the pseudo-random numbers' seed, etc., may lead to an entirely different evolution of the underlying search algorithms and hence, bring variability in the final performance. In \textit{AirCROP}, there exists two sources of variability: (i) the permutation of the newly generated pairings due to adoption of a parallel architecture, and (ii) the use of numerical seed for generation of pseudo-random numbers. While the former source of variability is eliminated by using a unique sorting criterion to sort the newly generated set of pairings, the latter source of variability is eliminated by using a fixed numerical seed for generating pseudo-random numbers. For more details on the sources of performance variability in \textit{AirCROP}, interested readers are referred to \citet{aggarwal2020aircrop}.
%
%
\par In the previous section, the results presented in Table~\ref{tab:allCGs} were restricted to the \textit{AirCROP}-runs against a fixed random number's seed equal to 0 (in there, the newly generated set of pairings were uniquely sorted, implying that the causes of performance variability in case of a repeated run were negated). Hence, it is intriguing to investigate the impact of a change in the random numbers' seed on the performance of the proposed CG heuristic with respect to CG Random and CG Standard. Table~\ref{tab:statResults} attempts to shed light on this question via empirical evidence for all the five test cases, wherein, \textit{for each test case}:
\begin{itemize}
    \item the \textit{AirCROP}'s performance variability owing to just the variation in the random number's seed is investigated. Towards it, five different seeds, denoted by \#, where, \# $\in$ \{0, 1, 2, 3, 4\} are used. Notably, a column marked by ``S$\rightarrow$\#'' presents the final CPOP/IPP solutions offered by  \textit{AirCROP}, against the random numbers' seed equal to \#. It may be also be noted that the scope of this investigation is restricted to five seeds under practical considerations of the experimental time (a single run may span upto 30 hours)
    %
    %
    \item the change in numerical seed does not seem to significantly affect the Cost-quality of the final \textit{AirCROP} solutions (for all three CG mechanisms). However, the variation in the associated run-time is stronger, and that could be attributed to different search trajectories corresponding to different permutations of variables or different random numbers.
    \item in general, both the mean ($\mu$) and standard deviation ($\sigma$) of the final Cost computed over the five seeds, are lower/better for the proposed CG heuristic compared to the alternatives in CG Random and CG Standard.
    %
\end{itemize}
The above established fact that the proposed CG heuristic can offer final CPOP/IPP solutions with comparable cost quality over the varying numerical seeds, and significantly better cost-quality than the other two CG mechanisms (CG Random \& CG Standard) endorses its robustness and effectiveness.
\begin{table}[htbp]
	\scriptsize
	\centering
	\caption{\textit{AirCROP}'s performance$^\ast$ with varying random numbers' seed: Proposed CG heuristic vis-$\grave{a}$-vis CG Random and CG Standard. A Column marked by ``S$\rightarrow$\#'' corresponds to the pseudo-random number seed \#. }
	\begin{tabular}{M{7mm}M{18mm}M{7mm}M{11mm}M{11mm}M{11mm}M{11mm}M{11mm}M{25mm}}
		\toprule
		\textbf{Test} & \multicolumn{2}{c}{\textbf{Underlying CG}} & \multicolumn{5}{c}{\textbf{Runs with different random numbers' seeds}} & \multirow{2.5}{*}{$\bm{\mu \pm \sigma}$} \\
		\cmidrule(lr){4-8}
		\textbf{Case} & \multicolumn{2}{c}{\textbf{mechanism}} & \textbf{S$\rightarrow$0} & \textbf{S$\rightarrow$1} & \textbf{S$\rightarrow$2} & \textbf{S$\rightarrow$3} & \textbf{S$\rightarrow$4} &\\
		\midrule
		\multirow{9}{*}{\textbf{TC-1}} & \multirow{2.5}{*}{\textbf{CG heuristic}} & \textbf{Cost} & 3471975 & 3470336 & 3472133 & 3470552 & 3468501 & 3470700 $\pm$ 1473 \\
		\cmidrule(lr){3-9} 
		&  & \textbf{Time} & 11:37 & 09:34 & 18:58 & 14:35 & 11:51 & 13:19 $\pm$ 03:38 \\
		\cmidrule(lr){2-9} 
		& \multirow{2.5}{*}{\textbf{CG Random}} & \textbf{Cost} & 3534145 & 3534649 & 3535918 & 3535263 & 3519618 & 3531919 $\pm$ 6909 \\
		\cmidrule(lr){3-9} 
		&  & \textbf{Time} & 16:41 & 16:50 & 18:02 & 17:03 & 19:49 & 17:41 $\pm$ 01:19 \\
		\cmidrule(lr){2-9} 
		& \multirow{2.5}{*}{\textbf{CG Standard}} & \textbf{Cost} & 3946632 & 3943849 & 3851649 & 3921841 & 3787803 & 3890355 $\pm$ 68984 \\
		\cmidrule(lr){3-9} 
		&  & \textbf{Time} & 07:36 & 07:01 & 07:02 & 05:32 & 07:15 & 06:54 $\pm$ 00:48 \\
		\midrule
		\multirow{9}{*}{\textbf{TC-2}} & \multirow{2.5}{*}{\textbf{CG heuristic}} & \textbf{Cost} & 3497010 & 3500929 & 3500478 & 3503829 & 3502424 & 3500934 $\pm$ 2560 \\
		\cmidrule(lr){3-9} 
		&  & \textbf{Time} & 15:22 & 12:24 & 13:11 & 11:17 & 12:37 & 12:59 $\pm$ 01:31 \\
		\cmidrule(lr){2-9} 
		& \multirow{2.5}{*}{\textbf{CG Random}} & \textbf{Cost} & 3536578 & 3547865 & 3544912 & 3546765 & 3546783 & 3544581 $\pm$ 4598 \\
		\cmidrule(lr){3-9} 
		&  & \textbf{Time} & 18:26 & 20:32 & 29:00 & 26:04 & 29:15 & 24:40 $\pm$ 04:57 \\
		\cmidrule(lr){2-9} 
		& \multirow{2.5}{*}{\textbf{CG Standard}} & \textbf{Cost} & 3902100 & 3926743 & 3915609 & 3814541 & 3892014 & 3890202 $\pm$ 44301 \\
		\cmidrule(lr){3-9} 
		&  & \textbf{Time} & 11:24 & 11:05 & 11:16 & 15:11 & 11:57 & 12:11 $\pm$ 01:43 \\
		\midrule
		\multirow{9}{*}{\textbf{TC-3}} & \multirow{2.5}{*}{\textbf{CG heuristic}} & \textbf{Cost} & 3495142 & 3487873 & 3488725 & 3488906 & 3486707 & 3489471 $\pm$ 3288 \\
		\cmidrule(lr){3-9} 
		&  & \textbf{Time} & 15:35 & 17:00 & 18:37 & 18:22 & 14:41 & 16:51 $\pm$ 01:43 \\
		\cmidrule(lr){2-9} 
		& \multirow{2.5}{*}{\textbf{CG Random}} & \textbf{Cost} & 3524914 & 3607717 & 3616467 & 3565324 & 3594284 & 3581742 $\pm$ 37205 \\
		\cmidrule(lr){3-9} 
		&  & \textbf{Time} & 28:28 & 32:44 & 30:56 & 32:38 & 27:51 & 30:32 $\pm$ 02:18 \\
		\cmidrule(lr){2-9} 
		& \multirow{2.5}{*}{\textbf{CG Standard}} & \textbf{Cost} & 3991375 & 3831094 & 3794025 & 3920094 & 3912296 & 3889777 $\pm$ 78040 \\
		\cmidrule(lr){3-9} 
		&  & \textbf{Time} & 10:49 & 17:14 & 19:36 & 10:51 & 08:58 & 13:30 $\pm$ 04:39 \\
		\midrule
		\multirow{9}{*}{\textbf{TC-4}} & \multirow{2.5}{*}{\textbf{CG heuristic}} & \textbf{Cost} & 3599539 & 3601667 & 3603315 & 3607525 & 3603422 & 3603094 $\pm$ 2936 \\
		\cmidrule(lr){3-9} 
		&  & \textbf{Time} & 12:41 & 12:46 & 10:59 & 11:00 & 10:08 & 11:31 $\pm$ 01:10 \\
		\cmidrule(lr){2-9} 
		& \multirow{2.5}{*}{\textbf{CG Random}} & \textbf{Cost} & 3643191 & 3686123 & 3656830 & 3651730 & 3646299 & 3656835 $\pm$ 17183 \\
		\cmidrule(lr){3-9} 
		&  & \textbf{Time} & 16:24 & 14:23 & 17:19 & 16:25 & 16:10 & 16:09 $\pm$ 01:05 \\
		\cmidrule(lr){2-9} 
		& \multirow{2.5}{*}{\textbf{CG Standard}} & \textbf{Cost} & 4093020 & 4006554 & 4036782 & 4032261 & 4020988 & 4037921 $\pm$ 32939 \\
		\cmidrule(lr){3-9} 
		&  & \textbf{Time} & 05:46 & 07:31 & 06:30 & 07:04 & 07:32 & 06:53 $\pm$ 00:46 \\
		\midrule
		\multirow{9}{*}{\textbf{TC-5}} & \multirow{2.5}{*}{\textbf{CG heuristic}} & \textbf{Cost} & 4593256 & 4595372 & 4593355 & 4592388 & 4594304 & 4593735 $\pm$ 1140 \\
		\cmidrule(lr){3-9} 
		&  & \textbf{Time} & 26:28 & 23:21 & 23:15 & 25:57 & 22:54 & 24:23 $\pm$ 01:42 \\
		\cmidrule(lr){2-9} 
		& \multirow{2.5}{*}{\textbf{CG Random}} & \textbf{Cost} & 4749477 & 4714886 & 4817132 & 4773801 & 4772414 & 4765542 $\pm$ 37421 \\
		\cmidrule(lr){3-9} 
		&  & \textbf{Time} & 23:10 & 21:07 & 17:14 & 22:02 & 18:00 & 20:19 $\pm$ 02:35 \\
		\cmidrule(lr){2-9} 
		& \multirow{2.5}{*}{\textbf{CG Standard}} & \textbf{Cost} & 5109643 & 5021205 & 5193059 & 5297016 & 5055532 & 5135291 $\pm$ 111274 \\
		\cmidrule(lr){3-9} 
		&  & \textbf{Time} & 31:38 & 32:44 & 20:20 & 22:49 & 34:12 & 28:21 $\pm$ 06:19 \\
		\bottomrule
	\end{tabular}
	\label{tab:statResults}
	\\\vspace{2pt}\footnotesize{$^\ast$All values in the ``Cost'' rows are in USD, and all the corresponding real values are rounded-off to the next integer values. All values in the ``Time'' rows are in HH:MM, and all the corresponding seconds' values are rounded-off to the next minute values. Notably, $\mu$ and $\sigma$ represent the mean and standard deviation, respectively, of the Cost or Time (over the five numerical seeds for random numbers), depending on the label of the row.}
\end{table}
\subsubsection{\textit{AirCROP} performance: Proposed CG heuristic vis-\texorpdfstring{$\grave{a}$}{a'}-vis Sub-combinations of its constituent CG strategies}
This section pertains to the second goal of investigating the performance of the proposed CG heuristic in contrast to a scenario where one or two of its constituent strategies may be done away with. To enable meaningful inferences, results for two test cases with varying flight-set sizes, namely, TC-2 (3228 flights) and TC-5 (4212 flights), are presented in Table~\ref{tab:CGtc2} and Table~\ref{tab:CGtc5}, respectively. Both these tables adopt an identical structure, in that:
\begin{itemize}
	\item the first row highlights the cost associated with the $IFS$ that triggers the \textit{AirCROP}.
	\item the last (emboldened) row marked by ``Final Solution'' highlights the results of the final CPOP/IPP solutions of the respective runs, in terms of cost in the last permissible interaction, and total run-time consumed.
	\item the intermediate rows facilitate a comparison between the results obtained from the proposed CG heuristic vis-$\grave{a}$-vis sub-combinations of its constituent CG strategies, in terms of the solution's cost quality and run-time consumed, across different LPP-solutioning and IPP-solutioning phases.
	\item notably, some sub-combinations of the four CG strategies require more LPP-IPP interactions than those which could be vertically accommodated. In such instances, given the paucity of space, only the first and last six LPP-IPP interactions are presented.
	\item for brevity, the CG strategies -- $CGD$, $CGU$, $CGA$ \& $CGR$, have been notated as $D,~U,~A~\&~R$, respectively. 
\end{itemize}
\afterpage{%
	\clearpage
	\thispagestyle{empty}
	\begin{landscape}
		\centering 
		\begin{table}[htbp]
			\scriptsize
			\caption{\textit{AirCROP} performance$^\ast$ on TC-2: Proposed CG heuristic vis-$\grave{a}$-vis Sub-combinations of its constituent CG strategies}
			\resizebox{\columnwidth}{!}{%
				\begin{tabular}{M{2mm}M{6mm}M{7mm}M{4mm}M{7mm}M{4mm}M{7mm}M{4mm}M{7mm}M{4mm}M{7mm}M{4mm}M{7mm}M{4mm}M{7mm}M{4mm}M{7mm}M{4mm}M{7mm}M{4mm}M{7mm}M{4mm}M{7mm}M{4mm}}
					\toprule
					\multicolumn{2}{c}{\multirow{1.25}{*}{\textbf{LPP-IPP}}} & \multicolumn{22}{c}{\textbf{$IFS$: 129221508 USD}} \\
					\cmidrule(l){3-24}
					\multicolumn{2}{c}{\multirow{1}{*}{\textbf{Interactions}}} & \multicolumn{2}{M{17mm}}{$D+U+A+R$} & \multicolumn{2}{M{17mm}}{$D+U$} & \multicolumn{2}{M{17mm}}{$D+A$} & \multicolumn{2}{M{17mm}}{$D+R$} & \multicolumn{2}{M{17mm}}{$U+A$} & \multicolumn{2}{M{17mm}}{$U+R$} & \multicolumn{2}{M{17mm}}{$A+R$} & \multicolumn{2}{M{17mm}}{$D+U+A$} & \multicolumn{2}{M{17mm}}{$D+U+R$} & \multicolumn{2}{M{17mm}}{$D+A+R$} & \multicolumn{2}{M{17mm}}{$U+A+R$} \\ \cmidrule(lr){1-2} \cmidrule(lr){3-4} \cmidrule(lr){5-6} \cmidrule(lr){7-8} \cmidrule(lr){9-10} \cmidrule(lr){11-12} \cmidrule(lr){13-14} \cmidrule(lr){15-16} \cmidrule(lr){17-18} \cmidrule(lr){19-20} \cmidrule(lr){21-22} \cmidrule(l){23-24}
					$T$ & $\mathcal{P}_{\#}^T$ & \textbf{Cost} & \textbf{Time} & \textbf{Cost} & \textbf{Time} & \textbf{Cost} & \textbf{Time} & \textbf{Cost} & \textbf{Time} & \textbf{Cost} & \textbf{Time} & \textbf{Cost} & \textbf{Time} & \textbf{Cost} & \textbf{Time} & \textbf{Cost} & \textbf{Time} & \textbf{Cost} & \textbf{Time} & \textbf{Cost} & \textbf{Time} & \textbf{Cost} & \textbf{Time} \\
					\midrule
					\multirow{2}{*}{$1$} & $\mathcal{P}_{LP}^1$ & 3495054 & 03:58 & 3579674 & 05:05 & 3599693 & 05:40 & 3503584 & 02:51 & 3734292 & 03:01 & 3506357 & 02:00 & 3504227 & 04:03 & 3554303 & 03:24 & 3495788 & 03:01 & 3501365 & 03:23 & 3497592 & 03:05 \\ \cmidrule{2-24}
					& $\mathcal{P}_{IP}^1$ & 3779111 & 00:20 & 3866200 & 00:20 & 3848648 & 00:20 & 3776618 & 00:20 & 4203521 & 00:20 & 3778414 & 00:20 & 3901568 & 00:20 & 3780826 & 00:20 & 3717868 & 00:20 & 3775231 & 00:20 & 3802960 & 00:20 \\ \midrule
					\multirow{2}{*}{$2$} & $\mathcal{P}_{LP}^2$ & 3492575 & 03:34 & 3641696 & 01:11 & 3604831 & 00:39 & 3514219 & 01:29 & 3688272 & 01:37 & 3509512 & 01:17 & 3504700 & 02:17 & 3556947 & 00:47 & 3502378 & 00:52 & 3501812 & 02:41 & 3497956 & 02:11 \\ \cmidrule{2-24}
					& $\mathcal{P}_{IP}^2$ & 3632511 & 00:20 & 3699026 & 00:09 & 3716080 & 00:20 & 3601799 & 00:20 & 3989067 & 00:20 & 3631382 & 00:20 & 3782327 & 00:20 & 3662331 & 00:20 & 3568175 & 00:20 & 3688671 & 00:20 & 3653415 & 00:20 \\ \midrule
					\multirow{2}{*}{$3$} & $\mathcal{P}_{LP}^3$ & 3494078 & 01:09 & 3625970 & 00:37 & 3600020 & 00:45 & 3580176 & 00:17 & 3675270 & 00:45 & 3523561 & 00:53 & 3503305 & 02:33 & 3561815 & 00:24 & 3503823 & 01:12 & 3500671 & 02:06 & 3498722 & 01:36 \\ \cmidrule{2-24}
					& $\mathcal{P}_{IP}^3$ & 3573900 & 00:20 & 3664508 & 00:01 & 3711086 & 00:20 & 3582853 & 00:01 & 3859771 & 00:20 & 3554587 & 00:01 & 3706886 & 00:20 & 3601047 & 00:01 & 3523485 & 00:01 & 3634846 & 00:20 & 3600274 & 00:20 \\ \midrule
					\multirow{2}{*}{$4$} & $\mathcal{P}_{LP}^4$ & 3494145 & 01:22 & 3622976 & 00:20 & 3599012 & 00:36 & 3579224 & 00:14 & 3655759 & 00:42 & 3525192 & 00:13 & 3503362 & 02:42 & 3560466 & 00:30 & 3505263 & 00:31 & 3500821 & 02:04 & 3500576 & 00:47 \\ \cmidrule{2-24}
					& $\mathcal{P}_{IP}^4$ & 3532658 & 00:01 & 3673371 & 00:02 & 3695672 & 00:20 & 3579224 & 00:01 & 3789943 & 00:20 & 3554187 & 00:01 & 3690647 & 00:20 & 3589096 & 00:01 & 3505263 & 00:01 & 3600917 & 00:20 & 3546205 & 00:02 \\ \midrule
					\multirow{2}{*}{$5$} & $\mathcal{P}_{LP}^5$ & 3495669 & 00:55 & 3630305 & 00:09 & 3597638 & 00:39 &  &  & 3658480 & 00:23 & 3525651 & 00:15 & 3503343 & 02:09 & 3560839 & 00:26 &  &  & 3501668 & 02:07 & 3501748 & 00:47 \\ \cmidrule{2-24}
					& $\mathcal{P}_{IP}^5$ & 3529639 & 00:01 & 3677469 & 00:02 & 3661911 & 00:13 &  &  & 3729362 & 00:20 & 3547749 & 00:01 & 3664605 & 00:20 & 3578777 & 00:01 &  &  & 3584076 & 00:20 & 3525674 & 00:01 \\ \midrule
					\multirow{2}{*}{$6$} & $\mathcal{P}_{LP}^6$ & 3496165 & 01:10 & 3633204 & 00:14 & 3597346 & 00:39 &  &  & 3653819 & 00:18 & 3523721 & 00:43 & 3504793 & 01:57 & 3558547 & 00:35 &  &  & 3501914 & 02:31 & 3502738 & 00:57 \\ \cmidrule{2-24}
					& $\mathcal{P}_{IP}^6$ & 3514831 & 00:01 & 3655049 & 00:01 & 3658258 & 00:10 &  &  & 3745261 & 00:20 & 3545140 & 00:01 & 3635246 & 00:20 & 3594226 & 00:01 &  &  & 3572761 & 00:20 & 3523248 & 00:01 \\ \midrule
					\multirow{2}{*}{$7$} & $\mathcal{P}_{LP}^7$ & 3497105 & 00:52 & \vdots & \vdots & 3596608 & 00:31 &  &  & \vdots & \vdots & 3523152 & 00:25 & 3503302 & 02:36 & 3560474 & 00:20 &  &  & 3504112 & 01:35 & 3503861 & 00:37 \\ \cmidrule{2-24} 
					& $\mathcal{P}_{IP}^7$ & 3501399 & 00:01 & \vdots & \vdots & 3647007 & 00:02 &  &  & \vdots & \vdots & 3550768 & 00:01 & 3633129 & 00:20 & 3572020 & 00:01 &  &  & 3527054 & 00:01 & 3512690 & 00:01 \\ \midrule
					&  &  &  & \multicolumn{2}{M{17mm}}{($T =$ 17 onwards $\downarrow$)} &  &  &  &  & \multicolumn{2}{M{17mm}}{($T =$ 19 onwards $\downarrow$)} &  &  &  &  &  &  &  &  &  &  &  &  \\ \midrule
					\multirow{2}{*}{$8$} & $\mathcal{P}_{LP}^8$ & 3496510 & 00:40 & 3640901 & 00:12 & 3601091 & 00:20 &  &  & 3633622 & 00:21 & 3524088 & 00:15 & 3503316 & 02:33 & 3558900 & 00:33 &  &  & 3514975 & 00:59 & 3505418 & 00:32 \\ \cmidrule{2-24}
					& $\mathcal{P}_{IP}^8$ & 3513764 & 00:01 & 3677841 & 00:01 & 3619447 & 00:01 &  &  & 3673195 & 00:01 & 3530800 & 00:01 & 3597729 & 00:20 & 3571622 & 00:01 &  &  & 3514975 & 00:01 & 3505418 & 00:01 \\ \midrule
					\multirow{2}{*}{$9$} & $\mathcal{P}_{LP}^9$ & 3497010 & 00:37 & 3646196 & 00:14 & 3604617 & 00:10 &  &  & 3632717 & 00:15 & 3523964 & 00:15 & 3504424 & 01:01 & 3559898 & 00:23 &  &  &  &  &  &  \\ \cmidrule{2-24}
					& $\mathcal{P}_{IP}^9$ & 3497010 & 00:01 & 3657922 & 00:01 & 3604617 & 00:01 &  &  & 3665437 & 00:01 & 3523964 & 00:01 & 3580000 & 00:20 & 3567262 & 00:01 &  &  &  &  &  &  \\ \midrule
					\multirow{2}{*}{$10$} & $\mathcal{P}_{LP}^{10}$ &  &  & 3646027 & 00:10 &  &  &  &  & 3632327 & 00:18 &  &  & 3503205 & 01:31 & 3560154 & 00:18 &  &  &  &  &  &  \\ \cmidrule{2-24}
					& $\mathcal{P}_{IP}^{10}$ &  &  & 3654784 & 00:01 &  &  &  &  & 3659656 & 00:01 &  &  & 3579938 & 00:20 & 3568940 & 00:01 &  &  &  &  &  &  \\ \midrule
					\multirow{2}{*}{$11$} & $\mathcal{P}_{LP}^{11}$ &  &  & 3646848 & 00:11 &  &  &  &  & 3633160 & 00:22 &  &  & 3505119 & 00:55 & 3561733 & 00:10 &  &  &  &  &  &  \\ \cmidrule{2-24}
					& $\mathcal{P}_{IP}^{11}$ &  &  & 3651184 & 00:01 &  &  &  &  & 3651975 & 00:01 &  &  & 3563265 & 00:20 & 3561733 & 00:01 &  &  &  &  &  &  \\ \midrule
					\multirow{2}{*}{$12$} & $\mathcal{P}_{LP}^{12}$ &  &  & 3643297 & 00:15 &  &  &  &  & 3635032 & 00:09 &  &  & 3504626 & 01:14 &  &  &  &  &  &  &  &  \\ \cmidrule{2-24}
					& $\mathcal{P}_{IP}^{12}$ &  &  & 3658888 & 00:01 &  &  &  &  & 3637733 & 00:01 &  &  & 3551695 & 00:20 &  &  &  &  &  &  &  &  \\ \midrule
					\multirow{2}{*}{$13$} & $\mathcal{P}_{LP}^{13}$ &  &  & 3645647 & 00:06 &  &  &  &  & 3635138 & 00:12 &  &  & 3528315 & 00:22 &  &  &  &  &  &  &  &  \\ \cmidrule{2-24}
					& $\mathcal{P}_{IP}^{13}$ &  &  & 3645647 & 00:01 &  &  &  &  & 3635138 & 00:01 &  &  & 3528315 & 00:01 &  &  &  &  &  &  &  &  \\ \midrule
					\multicolumn{2}{c}{\textbf{Final Solution}} & \textbf{3497010} & \textbf{15:23} & \textbf{3645647} & \textbf{11:32} & \textbf{3604617} & \textbf{11:41} & \textbf{3579224} & \textbf{05:33} & \textbf{3635138} & \textbf{16:02} & \textbf{3523964} & \textbf{07:03} & \textbf{3528315} & \textbf{29:54} & \textbf{3561733} & \textbf{08:39} & \textbf{3505263} & \textbf{06:18} & \textbf{3514975} & \textbf{19:28} & \textbf{3505418} & \textbf{11:38}\\
					\bottomrule
			\end{tabular}}
			\label{tab:CGtc2}
			\\\scriptsize{$^\ast$In this table: (a) for brevity, $D,~U,~A~\&~R$ are used to represent $CGD,~CGU,~CGA~\&~CGR$, respectively, (b) the unit for the ``Cost'' columns is USD (in that, the real values are rounded-off to the next integer), (c) the format for the ``Time'' columns is HH:MM (in that, the seconds' values are rounded-off to the next minute values), and (d) given the paucity of space, for some sub-combinations of strategies, only the first and last six LPP-IPP interactions are presented.}
		\end{table}
	\end{landscape}
	\clearpage
}
\afterpage{%
	\clearpage
	\thispagestyle{empty}
	\begin{landscape}
		\centering 
		\begin{table}[htbp]
			\scriptsize
			\caption{\textit{AirCROP} performance$^\ast$ on TC-5: Proposed CG heuristic vis-$\grave{a}$-vis Sub-combinations of its constituent CG strategies}
			\resizebox{\columnwidth}{!}{%
				\begin{tabular}{M{2mm}M{6mm}M{7mm}M{4mm}M{7mm}M{4mm}M{7mm}M{4mm}M{7mm}M{4mm}M{7mm}M{4mm}M{7mm}M{4mm}M{7mm}M{4mm}M{7mm}M{4mm}M{7mm}M{4mm}M{7mm}M{4mm}M{7mm}M{4mm}}
					\toprule
					\multicolumn{2}{c}{\multirow{1.25}{*}{\textbf{LPP-IPP}}} & \multicolumn{22}{c}{\textbf{$IFS$: 131443284 USD}} \\
					\cmidrule(l){3-24}
					\multicolumn{2}{c}{\multirow{1}{*}{\textbf{Interactions}}} & \multicolumn{2}{M{17mm}}{$D+U+A+R$} & \multicolumn{2}{M{17mm}}{$D+U$} & \multicolumn{2}{M{17mm}}{$D+A$} & \multicolumn{2}{M{17mm}}{$D+R$} & \multicolumn{2}{M{17mm}}{$U+A$} & \multicolumn{2}{M{17mm}}{$U+R$} & \multicolumn{2}{M{17mm}}{$A+R$} & \multicolumn{2}{M{17mm}}{$D+U+A$} & \multicolumn{2}{M{17mm}}{$D+U+R$} & \multicolumn{2}{M{17mm}}{$D+A+R$} & \multicolumn{2}{M{17mm}}{$U+A+R$} \\ \cmidrule(lr){1-2} \cmidrule(lr){3-4} \cmidrule(lr){5-6} \cmidrule(lr){7-8} \cmidrule(lr){9-10} \cmidrule(lr){11-12} \cmidrule(lr){13-14} \cmidrule(lr){15-16} \cmidrule(lr){17-18} \cmidrule(lr){19-20} \cmidrule(lr){21-22} \cmidrule(l){23-24}
					$T$ & $\mathcal{P}_{\#}^T$ & \textbf{Cost} & \textbf{Time} & \textbf{Cost} & \textbf{Time} & \textbf{Cost} & \textbf{Time} & \textbf{Cost} & \textbf{Time} & \textbf{Cost} & \textbf{Time} & \textbf{Cost} & \textbf{Time} & \textbf{Cost} & \textbf{Time} & \textbf{Cost} & \textbf{Time} & \textbf{Cost} & \textbf{Time} & \textbf{Cost} & \textbf{Time} & \textbf{Cost} & \textbf{Time} \\
					\midrule
					\multirow{2}{*}{$1$} & $\mathcal{P}_{LP}^1$ & 4584526 & 07:40 & 4676660 & 10:44 & 4757655 & 10:24 & 4620290 & 06:13 & 4995573 & 05:49 & 4596698 & 05:36 & 4608773 & 07:37 & 4652569 & 07:58 & 4589404 & 05:55 & 4593387 & 09:07 & 4591268 & 07:07 \\ \cmidrule(l){2-24}
					& $\mathcal{P}_{IP}^1$ & 4997918 & 00:20 & 5026401 & 00:20 & 5043653 & 00:20 & 5037950 & 00:20 & 5671588 & 00:20 & 5103737 & 00:20 & 5272271 & 00:20 & 5057001 & 00:20 & 4999238 & 00:20 & 4959504 & 00:20 & 5110753 & 00:20 \\ \midrule
					\multirow{2}{*}{$2$} & $\mathcal{P}_{LP}^2$ & 4590082 & 03:16 & 4751749 & 02:04 & 4758993 & 01:45 & 4639160 & 02:40 & 4938333 & 03:25 & 4610225 & 02:45 & 4601101 & 06:49 & 4651835 & 02:17 & 4603615 & 03:47 & 4595702 & 06:01 & 4590353 & 06:09 \\ \cmidrule(l){2-24}
					& $\mathcal{P}_{IP}^2$ & 4795955 & 00:20 & 4836224 & 00:20 & 4913084 & 00:20 & 4737834 & 00:20 & 5552726 & 00:20 & 4797107 & 00:20 & 5076692 & 00:20 & 4826001 & 00:20 & 4714892 & 00:20 & 4750037 & 00:20 & 4910403 & 00:20 \\ \midrule
					\multirow{2}{*}{$3$} & $\mathcal{P}_{LP}^3$ & 4589793 & 03:45 & 4761415 & 00:20 & 4766677 & 01:17 & 4714329 & 00:25 & 4903807 & 02:35 & 4639790 & 02:02 & 4602620 & 05:57 & 4653992 & 02:02 & 4626447 & 01:46 & 4596653 & 04:26 & 4592982 & 04:33 \\ \cmidrule(l){2-24}
					& $\mathcal{P}_{IP}^3$ & 4743617 & 00:20 & 4794472 & 00:01 & 4831406 & 00:06 & 4714329 & 00:01 & 5379158 & 00:20 & 4697779 & 00:01 & 5030872 & 00:20 & 4739588 & 00:20 & 4653897 & 00:01 & 4676496 & 00:20 & 4765292 & 00:20 \\ \midrule
					\multirow{2}{*}{$4$} & $\mathcal{P}_{LP}^4$ & 4592119 & 02:41 & 4758381 & 00:26 & 4793764 & 00:34 &  &  & 4878453 & 02:01 & 4641037 & 01:44 & 4603899 & 06:09 & 4658183 & 00:46 & 4621255 & 00:33 & 4596886 & 04:19 & 4594583 & 02:44 \\ \cmidrule(l){2-24}
					& $\mathcal{P}_{IP}^4$ & 4685270 & 00:20 & 4800989 & 00:01 & 4793764 & 00:01 &  &  & 5211053 & 00:20 & 4698508 & 00:13 & 4992284 & 00:20 & 4689930 & 00:01 & 4631162 & 00:01 & 4642862 & 00:02 & 4663171 & 00:20 \\ \midrule
					\multirow{2}{*}{$5$} & $\mathcal{P}_{LP}^5$ & 4591363 & 02:19 & 4756242 & 00:21 &  &  &  &  & 4852866 & 01:48 & 4641852 & 01:10 &  &  & 4658672 & 00:33 & 4620508 & 00:36 & 4618622 & 01:04 & 4596460 & 01:59 \\ \cmidrule(l){2-24}
					& $\mathcal{P}_{IP}^5$ & 4672393 & 00:20 & 4792818 & 00:01 &  &  &  &  & 5103427 & 00:20 & 4663457 & 00:01 &  &  & 4687065 & 00:01 & 4631519 & 00:01 & 4618622 & 00:01 & 4630592 & 00:01 \\ \midrule
					\multirow{2}{*}{$6$} & $\mathcal{P}_{LP}^6$ & 4591374 & 02:26 & 4761148 & 00:21 &  &  &  &  & 4851587 & 01:03 & 4648513 & 00:29 &  &  & 4660279 & 00:29 & 4619001 & 00:56 &  &  & 4597909 & 01:50 \\ \cmidrule(l){2-24}
					& $\mathcal{P}_{IP}^6$ & 4631706 & 00:01 & 4784897 & 00:01 &  &  &  &  & 5017600 & 00:20 & 4668904 & 00:01 &  &  & 4672588 & 00:01 & 4662053 & 00:01 &  &  & 4630469 & 00:01 \\ \midrule
					\multirow{2}{*}{$7$} & $\mathcal{P}_{LP}^7$ & 4593502 & 01:28 & \vdots & \vdots &  &  &  &  & \vdots & \vdots & 4647303 & 00:50 &  &  & 4660895 & 00:21 & \vdots & \vdots &  &  & 4598584 & 01:49 \\ \cmidrule{2-24} 
					& $\mathcal{P}_{IP}^7$ & 4604572 & 00:01 & \vdots & \vdots &  &  &  &  & \vdots & \vdots & 4676948 & 00:01 &  &  & 4670384 & 00:01 & \vdots & \vdots &  &  & 4616375 & 00:01 \\ \midrule
					&  &  &  & \multicolumn{2}{M{17mm}}{($T =$ 25 onwards $\downarrow$)} &  &  &  &  & \multicolumn{2}{M{17mm}}{($T =$ 16 onwards $\downarrow$)} &  &  &  &  &  &  & \multicolumn{2}{M{17mm}}{($T =$ 12 onwards $\downarrow$)} &  &  &  &  \\ \midrule
					\multirow{2}{*}{$8$} & $\mathcal{P}_{LP}^8$ & 4593256 & 01:06 & 4755574 & 00:13 &  &  &  &  & 4840529 & 00:30 & 4650352 & 00:24 &  &  & 4660975 & 00:16 & 4624007 & 00:28 &  &  & 4599165 & 01:12 \\ \cmidrule(l){2-24}
					& $\mathcal{P}_{IP}^8$ & 4593256 & 00:01 & 4763690 & 00:01 &  &  &  &  & 4884804 & 00:02 & 4658859 & 00:01 &  &  & 4666961 & 00:01 & 4641306 & 00:01 &  &  & 4614854 & 00:01 \\ \midrule
					\multirow{2}{*}{$9$} & $\mathcal{P}_{LP}^9$ &  &  & 4748472 & 00:26 &  &  &  &  & 4840343 & 00:32 & 4646354 & 00:30 &  &  & 4661466 & 00:26 & 4624237 & 00:44 &  &  & 4600468 & 00:48 \\ \cmidrule(l){2-24}
					& $\mathcal{P}_{IP}^9$ &  &  & 4781479 & 00:01 &  &  &  &  & 4879210 & 00:01 & 4683759 & 00:01 &  &  & 4661466 & 00:01 & 4639496 & 00:01 &  &  & 4603901 & 00:01 \\ \midrule
					\multirow{2}{*}{$10$} & $\mathcal{P}_{LP}^{10}$ &  &  & 4748601 & 00:18 &  &  &  &  & 4835805 & 00:37 & 4652942 & 00:48 &  &  &  &  & 4626408 & 00:20 &  &  &  &  \\ \cmidrule(l){2-24}
					& $\mathcal{P}_{IP}^{10}$ &  &  & 4782221 & 00:01 &  &  &  &  & 4882052 & 00:02 & 4669279 & 00:01 &  &  &  &  & 4651319 & 00:01 &  &  &  &  \\ \midrule
					\multirow{2}{*}{$11$} & $\mathcal{P}_{LP}^{11}$ &  &  & 4753802 & 00:20 &  &  &  &  & 4835077 & 00:23 & 4654569 & 00:13 &  &  &  &  & 4624096 & 00:20 &  &  &  &  \\ \cmidrule(l){2-24}
					& $\mathcal{P}_{IP}^{11}$ &  &  & 4771304 & 00:01 &  &  &  &  & 4881412 & 00:03 & 4660035 & 00:01 &  &  &  &  & 4626663 & 00:01 &  &  &  &  \\ \midrule
					\multirow{2}{*}{$12$} & $\mathcal{P}_{LP}^{12}$ &  &  & 4754076 & 00:25 &  &  &  &  & 4835290 & 00:51 & 4651333 & 00:39 &  &  &  &  & 4624002 & 00:22 &  &  &  &  \\ \cmidrule(l){2-24}
					& $\mathcal{P}_{IP}^{12}$ &  &  & 4777057 & 00:01 &  &  &  &  & 4867568 & 00:01 & 4660973 & 00:01 &  &  &  &  & 4627987 & 00:01 &  &  &  &  \\ \midrule
					\multirow{2}{*}{$13$} & $\mathcal{P}_{LP}^{13}$ &  &  & 4753421 & 00:11 &  &  &  &  & 4835114 & 00:40 & 4653513 & 00:14 &  &  &  &  & 4623623 & 00:27 &  &  &  &  \\ \cmidrule(l){2-24}
					& $\mathcal{P}_{IP}^{13}$ &  &  & 4762025 & 00:01 &  &  &  &  & 4861977 & 00:01 & 4653513 & 00:01 &  &  &  &  & 4623623 & 00:01 &  &  &  &  \\ \midrule
					\multicolumn{2}{c}{\textbf{Final solution}} & \textbf{4593256} & \textbf{26:24} & \textbf{4762025} & \textbf{26:05} & \textbf{4793764} & \textbf{14:47} & \textbf{4714329} & \textbf{09:59} & \textbf{4861977} & \textbf{29:43} & \textbf{4653513} & \textbf{18:27} & \textbf{4992284} & \textbf{27:52} & \textbf{4661466} & \textbf{16:14} & \textbf{4623623} & \textbf{19:51} & \textbf{4618622} & \textbf{26:00} & \textbf{4603901} & \textbf{29:36} \\
					\bottomrule
			\end{tabular}}
			\label{tab:CGtc5}
			\\\scriptsize{$^\ast$In this table: (a) for brevity, $D,~U,~A~\&~R$ are used to represent $CGD,~CGU,~CGA~\&~CGR$, respectively, (b) the unit for the ``Cost'' columns is USD (in that, the real values are rounded-off to the next integer), (c) the format for the ``Time'' columns is HH:MM (in that, the seconds' values are rounded-off to the next minute values), and (d) given the paucity of space, for some sub-combinations of strategies, only the first and last six LPP-IPP interactions are presented.}
		\end{table}
	\end{landscape}
	\clearpage
}
Interestingly, both the Tables~\ref{tab:CGtc2}~\&~\ref{tab:CGtc5} reveal similar trends, reinforcing confidence in their robustness. The prominent common trends are, as highlighted below:  
\begin{enumerate}
	\item the final CPOP/IPP solution offered by the proposed CG heuristic is better in terms of the cost quality, compared to any sub-combination of the CG strategies involved (taken two or three at a time), though the required run-times are comparable. Notably, the cost difference between the best and the second-best solution instances is to the tune of 8,253 USD for TC-2, and 10,645 USD for TC-5. Manifestation of a similar cost benefit, despite a scaled-up flight-set size (in TC-5), in a sense, reflects on to the promise of scalability in the proposed CG heuristic.
	\item the final CPOP/IPP solution offered by any sub-combination of three CG strategies necessarily involving $CGR$ as one of the constituents, is better in cost, than every solution offered by any sub-combination of two CG strategies. This trend asserts the importance of incorporating \textit{random-search} as a constituent strategy. 
	%
	%
	%
	\item within the sub-combinations of three CG strategies: those which involve $CGR$ as one of the constituent strategy, offer better cost compared to those which do not involve $CGR$. This trend again asserts the importance of incorporating \textit{random search} as a constituent strategy.
	\item within the sub-combinations of two CG strategies: barring an exception of one instance (in TC-5), those which involve $CGR$ as one of the constituent strategy, offer better cost compared to those which do not involve $CGR$. Despite the singular exception, the importance of incorporating \textit{random search} as a constituent strategy cannot be discounted.  
\end{enumerate}

Based on the above, it is fair to infer that among all the CG strategies considered one at a time, $CGR$ stands out as the most important strategy. However, $CGR$ accompanied by exploitation of the optimal solution features through $CGD$, $CGU$ and partly $CGA$ (since it fosters $CGD$ and $CGU$ based pairings from previous iterations) manifests as the best combination, that we refer to as the proposed CG heuristic. Notably, the CG heuristic based solutions have been validated by the industrial sponsor, as surpassing the best-known for these flight sets. The main limitation linked to these results is the significant run time. However, this could be attributed to the use of Python scripting language, and it is fair to expect significant reduction in run time through use of other programming languages, such as C++.

\section{Conclusion and Future Research} \label{sec:conc}
Airline crew pairing optimization problem is perceived as one of the most challenging combinatorial optimization problems in the OR-domain. Since, crew operating cost is the second-largest operating cost for an airline, it is critically important to generate a set of legal crew pairings, covering the given flight schedule at minimum cost possible. Numerous contributions have been made by researchers in the past to tackle the flight-networks with evolving scale-and-complexity by developing solutions around OR-techniques. Despite this progress, the much prevalent and emergent \textit{complex flight-networks} largely remain uninvestigated. This is all the more alarming, considering that the air traffic is expected to grow double in 20 years \cite{marisa2018airtravel}, wherein, majority of airlines may incorporate multiple crew bases and multiple hubs in their flight-networks. This research has proposed a novel domain-knowledge driven CG heuristic for efficiently tackling real-world flight networks characterized by an unprecedented (conjunct) scale-and-complexity. This CG heuristic constitutes the core search-mechanism of an in-house developed optimizer, namely, \textit{AirCROP}, which has been tested and validated for such real-world flight networks. To solve a given CPOP, \textit{AirCROP} relies on intermittent interaction of two phases -- \textit{LPP-solutioning} (relies on relaxing the integer constraints, and solves the resulting LPP, to fetch a low cost LPP solution), and \textit{IPP-solutioning} (integerizes the resulting LPP solution). In that, \textit{the quality of the final IPP solution largely depends on the search-efficiency in the LPP-solutioning phase, at the core of which lies the proposed CG heuristic, bearing testament to the importance of this research.}
%
\par Beginning with an initial feasible solution, the proposed CG heuristic iteratively introduces pairings, which gradually enable improvement in the the cost quality of the solution. In this endeavor, it relies on balancing the \textit{random exploration} of the pairings' space ($CGR$) and exploitation of optimal solution features including --  \textit{minimal deadheads at a set level} ($CGD$) \& \textit{high crew utilization at a pairing level} ($CGU$), while utilizing the past computational effort guided by \textit{flight-pair level} information ($CGA$). The efficacy of the proposed CG heuristic has been investigated on five real-world, large-scale (over 4,212 flights), complex flight networks (over 15 crew bases and multiple hub-and-spoke sub-networks). Through experimental evidence, it is established that the proposed CG heuristic infuses better search efficiency compared to the stand-alone \textit{random-search}, stand-alone \textit{exact-search}, or any sub-combinations of these CG strategies. 
\par Aligned with the industrial sponsor's larger research scheme, the codes of \textit{AirCROP}, including the proposed CG heuristic, have been developed using Python scripting language. However, a significant reduction in run-time is expected via the use of other programming languages such as C++, etc. Moreover, certain thresholds that are used to guide the \textit{search} have been set by factoring the flight-network characteristics on the one hand, and the limitations of available computational resources. Development of methods or strategies that can help reduce the dependency of the final performance on threshold settings, may be an important future research direction. Furthermore, developing CG strategies during the IPP-solutioning may also augment the overall \textit{search} efficiency. Lastly, the emergent trend on utilizing the machine learning capabilities to assist combinatorial optimization tasks, may also hold promise for the airline crew pairing optimization. Despite the scope for improvement, the authors believe that this research could serve as a template on how to utilize the domain knowledge to enhance the \textit{search-efficiency} towards solving large-scale combinatorial optimization problems.
%

%
%
\section*{Acknowledgment}
This research work is an outcome of an Indo-Dutch joint research project, supported by the Ministry of Electronics and Information Technology (MEITY), India [grant 13(4)/2015-CC\&BT]; Netherlands Organization for Scientific Research (NWO), the Netherlands; and GE Aviation, India. The authors acknowledge the invaluable support of GE Aviation team members -- Arioli Arumugam (Senior Director- Data \& Analytics), and Alla Rajesh (Senior Staff Data \& Analytics Scientist), for providing the problem definition, real-world test cases, and sharing the domain-knowledge during several insightful discussions which have helped the authors in successfully completing this research work.
%
%
\bibliographystyle{abbrvnat}       
\bibliography{CG_heuristic_ref}   
\end{document}